\documentclass[11pt,reqno]{amsart}
\usepackage{amsmath, amsfonts, amssymb, amscd}

\usepackage{amsmath, amsfonts, amssymb, amscd, enumerate}
\theoremstyle{plain}
\newtheorem{theo}{Theorem}[section]
\newtheorem{lem}[theo]{Lemma}
\newtheorem{prop}[theo]{Proposition}

\newtheorem{cor}[theo]{Corollary}
\theoremstyle{definition}
\newtheorem{rem}[theo]{Remark}

\newtheorem{definition}[theo]{Definition}
\newenvironment{pf}{\noindent{\it Proof. }}{$\square$\par\medskip}
\newenvironment{pfns}{\noindent{\it Proof. }}{\par\medskip}

\theoremstyle{plain}

\theoremstyle{definition}

\renewcommand{\=}{:=}

\newcommand{\beq}{\begin{equation}}
\newcommand{\eeq}{\end{equation}}
\renewcommand{\a}{\alpha}

\renewcommand{\d}{\delta}

\newcommand{\g}{\gamma}

\renewcommand{\k}{\kappa}
\renewcommand{\l}{\lambda}
\renewcommand{\o}{\omega}
\newcommand{\q}{\vartheta}

\newcommand{\s}{\sigma}
\renewcommand{\t}{\tau}

\newcommand{\z}{\zeta}
\newcommand{\D}{\Delta}

\newcommand{\G}{\Gamma}

\renewcommand{\L}{\Lambda}
\renewcommand{\O}{\Omega}

\newcommand{\bB}{\mathbb{B}}
\newcommand{\bC}{\mathbb{C}}

\newcommand{\bR}{\mathbb{R}}

\newcommand{\cA}{\mathcal{A}}
\newcommand{\cB}{\mathcal{B}}
\newcommand{\cC}{\mathcal{C}}

\newcommand{\cF}{\mathcal{F}}
\newcommand{\cG}{\mathcal{G}}
\newcommand{\cH}{\mathcal{H}}

\newcommand{\cK}{\mathcal{K}}

\newcommand{\cU}{\mathcal{U}}
\newcommand{\cV}{\mathcal{V}}

\newcommand{\cZ}{\mathcal{Z}}

\newcommand{\Jst}{J_o}

\renewcommand{\square}{\kern1pt\vbox
{\hrule height 0.6pt\hbox{\vrule width 0.6pt\hskip 3pt
\vbox{\vskip 6pt}\hskip 3pt\vrule width 0.6pt}\hrule height0.6pt}\kern1pt}

\DeclareMathOperator\Id{Id}

\renewcommand\Re{\operatorname{Re}}
\renewcommand\Im{\operatorname{Im}}

\newcommand{\wt}{\widetilde}
\newcommand{\wh}{\widehat}

\newcommand{\n}{\nabla}

\newcommand{\be}{\begin{equation}}
\newcommand{\ee}{\end{equation}}

\def\<#1,#2>{\langle\,#1,\,#2\,\rangle}
\newcommand{\arr}{\begin{array}{rlll}}
\newcommand{\ea}{\end{array}}
\newcommand{\bea}{\begin{eqnarray}}
\newcommand{\eea}{\end{eqnarray}}
\newcommand{\bean}{\begin{eqnarray*}}
\newcommand{\eean}{\end{eqnarray*}}

\def\sideremark#1{\ifvmode\leavevmode\fi\vadjust{
\vbox to0pt{\hbox to 0pt{\hskip\hsize\hskip1em
\vbox{\hsize3cm\tiny\raggedright\pretolerance10000
\noindent #1\hfill}\hss}\vbox to8pt{\vfil}\vss}}}

\newcounter{ssig}
\newcounter{ttig}
\title[Regular Monge-Amp\`ere exhaustions]
{Modular data and  \\ regularity of Monge-Amp\`ere exhaustions \\ and of Kobayashi distance}
\author[G. Patrizio and A.  Spiro]{Giorgio Patrizio
and Andrea Spiro}


\thanks{{\it Acknowledgments}. This research was partially supported by the Project MIUR ``Real and Complex Manifolds: Geometry, Topology and  Harmonic Analysis'' and by GNSAGA of INdAM}

\begin{document}

\begin{abstract}   Regularity properties of intrinsic objects for a large class of Stein Manifolds, namely of Monge-Amp\`ere exhaustions and  Kobayashi distance, is  interpreted in terms of modular data. The results lead to a  construction of  an infinite dimensional family of convex domains   with squared Kobayashi distance of prescribed regularity properties. A new sharp refinement of Stoll's characterization of $\bC^{n}$ is also given.

\end{abstract}

\maketitle

\null
\setcounter{section}{0}
\section{Introduction}
A well known theorem of Stoll (\cite{St}) characterizes 
 the unit ball $\bB^n \subset \bC^n$   as the unique  (up to biholomorphisms) $n$-dimensional complex manifold, which admits  a  $C^{\infty}$ strictly plurisubharmonic exhaustion $\tau:\bB^n\longrightarrow  [0,1)$   such that  the function $u\=\log \tau|_{\{\t > 0\}}$ satisfies the complex homogeneous Monge-Amp\`ere equation 
$\left(dd^{c}u\right)^{n}=0$ on $\bB^n \setminus \{\t = 0\}$. \par
\smallskip
In this theorem, it is a crucial assumption the smoothness   of $\t$ at all points and, in particular  on the minimal  set $\{\tau=0\}$. In fact,  Lempert (\cite{Le})   proved that  for any smooth strictly linearly convex domain  $D \subset \bC^n$ and  any  $p\in D$,  there exists  a   strictly plurisubharmonic exhaustion $\tau^{(p)}:D\longrightarrow [0,1)$ 
such that   $u\=\log \tau^{(p)}|_{\{\t > 0\}}$ satisfies the Monge-Amp\`ere equation  $\left(dd^{c}u\right)^{n}=0$ and is  
 of class $C^{\infty}$ at  all points of  $D \setminus \{p\}$, but,  in general, not at  $ \{p\} = \{\t = 0\}$. \par
 \smallskip
In other words, from a rigid situation (the  manifold is necessarily  $\bB^n$), the loss of smoothness on the minimal set leads  to a quite dispersed one, 
represented by  the   infinite dimensional space of the  biholomorphic classes of the  strictly linearly convex domains of $\bC^n$. Indeed, the property that the biholomorphic classes of convex domains
constitute an  infinite dimensional space   was    precisely assessed by Lempert  in \cite{Le2} and  Bland and Duchamp in \cite{BD, BD1},   where    complete sets of modular data,   parameterizing  the   biholomorphic classes of pointed strictly  linearly convex domains,  were determined. \par
\smallskip
The  construction of Bland and Duchamp's modular data  was later simplified and extended in \cite{PS},  where  we proved that  they  actually   provide   a full moduli space for a large class of complex manifolds,  the so-called {\it manifolds of circular type}. This  class   was originally  introduced by the first author in 
\cite{Pt1} and properly includes all strictly linearly convex domains and all
 strictly pseudoconvex complete circular domains domains in $\bC^n$. Such manifolds are characterized by  the existence of  a Monge-Amp\`ere exhaustion, 
i.e.  a plurisubharmonic  exhaustion $\t$ which is smooth at all points  except possibly  at the minimal set $  \{\t = 0\}$ and   with  $u \= \log \t$  
satisfying  the Monge-Amp\`ere equation  (see \S 2.2 below, for the detailed definition).
\par
\smallskip
The existence of such Monge-Amp\`ere exhaustion  implies the existence of a very well behaved Green pluricomplex potential and it  is   tied with   another fundamental biholomorphic invariant,  the Kobayashi  distance. For instance,  the    exhaustion $\tau^{(p)}:D\longrightarrow  [0,1)$ of a strictly linearly convex domain $D$ is equal to the square of the  hyperbolic tangent  of the   Kobayashi distance from the point $p$ and  its second order truncation   coincides with the  squared Minkowski 
functional of the Kobayashi  indicatrix at $p$ (see e.g.  \cite{Pt3}). \par
\smallskip
These properties carry over to all manifolds of circular type and allow to relate the regularity  of  Monge-Amp\`ere exhaustions of these manifolds  with the regularity of the squared Kobayashi distance. 
\par
\smallskip
In  case of the unit ball $\bB^n$, the squared Kobayashi distance  from a point $p$ coincides with the squared hyperbolic distance determined by   the standard K\"ahler metric of constant negative curvature.   It is real analytic and the  Kobayashi indicatrix  at $p$ is the unit ball itself (up  to complex linear transformations). It is therefore natural to ask whether regularity properties of  Kobayashi distance, of  Kobayashi indicatrix  or  other intrinsic objects may be used to characterize certain classes of complex manifolds, such as the unit ball. \par
\smallskip
On the other hand  since long  it is  known that there exist strictly convex domains, which are not biholomorphic to $\bB^n$ and yet with the unit ball as indicatrix at some point.  This  was first discovered by  Sibony \cite{Si}  and  further analyzed  using  various kinds of modular data  in \cite{Le2}, \cite{BD}, \cite{BD1}, \cite{PS} (see also \S \ref{52} later).  There is  indeed   an infinite dimensional family of non-biholomorphic complex manifolds, all  with the unit ball as Kobayashi indicatrix at some point.
 \par
 \smallskip
In this paper we  consider the  problem of characterization of  the  manifolds of circular type (and,  consequently,  of circular domains and strictly linearly convex domains) admitting  Monge-Amp\`ere exhaustions and  Kobayashi distances of prescribed regularity. Our main tool  is the existence  of normal forms for manifolds of circular type, as proved in \cite{PS}.  They are  pairs $(\bB^n, J)$, formed by the standard   unit ball $\bB^n$ equipped with  a non-standard integrable almost complex structure $J$,   smoothly deformable  into the standard complex structure $\Jst$ and satisfying  appropriate additional conditions (see \S 2.3). In \cite{PS} it is shown that any manifold of circular type is biholomorphic to a normal form $(\bB^n, J)$ and any normal form   is in turn uniquely determined by a complex tensor field  $Ê\phi_J$ on the blow up of $\bB^n$ at the origin, which allows to express $J$ as a deformation of $\Jst$.  Using the properties of manifold of circular type, one can show that the deformation tensor $\phi_J$  can be uniquely decomposed as a series of the form $\phi_J   =  \sum_{k \geq 0}Ê\phi^{(k)}_J$, 
where the terms  $\phi_J^{(k)}$  are tensor fields satisfying appropriate  conditions and 
which we refer as {\it Bland-Duchamp invariants}. \par
\smallskip
It was remarked in \cite{PS} that the vanishing of all invariants $\phi_J^{(k)}$ with  $k\geq 1$ is equivalent to the fact that the considered manifold is (up to a biholomorphism)  a complete circular domain and coincide with its Kobayashi indicatrix at the center.  In this paper,  we  analyze  the somehow opposite situation, namely the vanishing   of an initial segment of   invariants, i.e.  $\phi^{(k)} = 0$  for all  $0 \leq k \leq k_o$. 
Our main results may be summarized as follows: 
\par
\vskip0.3cm
\begin{itemize}
\item[--] {\it if  $\phi^{(j)} = 0$ for all $0 \leq j \leq  k-1$ for some $k \geq 6$, then the Monge-Amp\`ere exhaustion $\t$ is of class $\cC^{\left[\frac{k}{2}\right]}$ at the center} (Theorem \ref{prop41}); 
\\
\item[--] {\it conversely, if the Monge-Amp\`ere exhaustion $\t$ is of class $\cC^{2m}$   at the center for some $m \geq 1$,  then  $\phi^{(j)} = 0$ for all $0 \leq j \leq  m -1$} (Theorem \ref{stoll}).
\end{itemize}

As a consequence,  the vanishing of the first six invariants $\phi^{(0)}_J$,  \ldots, $\phi^{(5)}_J$
implies that the Monge-Amp\`ere exhaustion
$\t$ is at least  $\cC^3$ at the center. From this follows  that the Kobayashi indicatrix
at the center is the unit ball and gives a new  precise 
quantitative way to see that,  even if one restricts   to  the special class of  strictly linearly convex domains in $\bC^n$, there exists an infinite dimensional family of domains, which  have the ball as Kobayashi indicatrix at a point, but are not  biholomorphic to the unit ball. 
Indeed, we are  able to show that,  for any $m \geq 3$, there exists a strictly linearly convex domain $D \subset \bC^n$ and a point $x_o \in D$, whose  associated  Monge-Amp\`ere exhaustion  is of  class  $\cC^{m}$ but not  $\cC^\infty$   (more precisely, not  $\cC^{4m + 2}$).  None of  such domains is biholomorphic to $\bB^n$, even though one can construct examples  that are arbitrarily  close to  $\bB^n$. 
\par
\smallskip
For strictly linearly convex domains, the Monge-Amp\`ere exhaustion $\t^{(p)}$, with center $p$,  is the square of the hyperbolic tangent of the Kobayashi  distance from $p$.  Hence  the above result implies that {\it for any  $m \geq 3$ there exist arbitrarily small deformations of $\bB^n$, for which  the squared Kobayashi distance from a point is of class $\cC^{m}$ but not   $\cC^{4m + 2}$}.
\par
\smallskip
We recall that, in \cite{Bu},   Burns presented an alternative  proof of 
 Stoll's characterization of $\bC^n$ in terms of (unbounded) Monge-Amp\`ere exhaustions, in which the    original smoothness assumptions are lowered  to $\cC^5$. Since then,  it looked  reasonable to   expect  a weakening  of  the regularity assumptions in  Stoll's characterization of  $\bB^n$.
The above examples show that  this is not 
 possible and that   in Stoll's characterization of  $\bB^n$,
the regularity assumptions   cannot be relaxed. \par
 \smallskip
 On the other hand,  our results  bring also to  a radical simplification of  the proof of  Stoll's characterization of $\bC^n$, 
 in which we may  
 further relax Burn's hypothesis and proving that 
 it  suffices that  the exhaustion is of class   $\cC^4$. \par
\medskip
The structure of the paper is the following. In  \S 2, we briefly recall  definitions and basic facts  on 
integrable almost complex structures, Monge-Amp\`ere exhaustions, manifolds of circular type and their normal forms.  In   \S 3, we  
clarify some delicate points on  complex manifold structures of normal forms,  which are crucial for the  later discussions. In \S 4 and \S 5, the main results are proved. The final  two sections are dedicated to the construction of   Monge-Amp\`ere exhaustions with prescribed
regularity at the center and to a new proof of Stoll's characterization of $\bC^n$, respectively. \par
\bigskip
{\it Notation.} Throughout the paper, $\bB^n$ is the unit ball of $\bC^n$ centered at the origin and $\D$ is the unit disk in $\bC$. We also denote 
by 
$e^o_1 = (1, 0\ldots, 0)$, \ldots, $e^o_n = (0, 0\ldots, 1)$
 the standard  basis of $\bC^n$ and by $\Jst$  the  standard complex structure of $\bC^n$, i.e. the  tensor field of type $(1,1)$ on $\bR^{2n} = \bC^n$ such that
$\Jst\left(\frac{\partial}{\partial x^i}\right) = \frac{\partial}{\partial y^i}$ and $\Jst\left(\frac{\partial}{\partial y^i}\right) = -\frac{\partial}{\partial x^i}$, where  $x^i = \Re(z^i)$ and $y^i = \Im(z^i)$.
\par

\bigskip
\section{Preliminaries}
\setcounter{equation}{0}
\subsection{Integrable almost complex structures and complex manifolds}\hfill\par
This section is devoted to recall very standard definitions and facts on complex structures. It aims to set  notations and to underline some fine points about the differentiability requirements, which are  crucial in what follows. 
\par
\smallskip
A {\it complex manifold $M$ of  dimension $n$}  is a topological manifold endowed with a complete  atlas $\cA$ of $\bC^n$-valued  homeomorphisms $\xi: \cU \subset M \longrightarrow \cV \subset \bC^n$ with  holomorphic overlaps. 
Newlander-Niremberg  Theorem establishes a natural correspondence between  complex manifold structures  and  formally integrable almost complex structures.  We briefly review this  important fact. \par
\medskip

An almost complex structure $J$   of class $\cC^{k+\a}$, $0 \leq \a < 1$, on a $2n$-dimensional real manifold $M$, is a family of endomorphisms
$$J_x: T_x M \longrightarrow T_x M \ , \ x \in M\ , \qquad \text{with}\qquad  J_x^2 = - \Id_{T_x M}\ ,$$
such that, denoting  $J_x =J^i_j\left.\frac{\partial}{\partial x^i}\otimes dx^j\right|_x$ for a given set of coordinates $(x^1, \ldots, x^{2n})$,  the components ${J^i_j}$ are   real  functions  of class $\cC^{k + \a}$. 
 If  $k\geq 1$,  the  {\it Nijenhuis tensor} $N_J$ is a tensor field  of type $(1,2)$,  whose value on   any  pair of smooth vector fields $X$, $Y$ is  
  $$N_J(X, Y) = [X, Y] - [JX, JY] + J[X, JY] - J[JX, Y]\ .$$
The almost complex structure $J$ is  called {\it formally integrable} if $N_J \equiv 0$. \par
\medskip
Given an almost complex structure,  each complexified tangent space $T^\bC_x M$, $x \in M$,  splits into  a direct sum 
$$T^\bC_x M = T^{1,0}_x M\oplus T^{0,1}_x M\ ,$$
where $ T^{1,0}_x M$ and  $T^{0,1}_x M = \overline{T^{1,0}_x M}$ are   the $(+i)$- and $(-i)$-eigenspace of $J_x$, respectively.  Around any point $x_o \in M$, there always exists  a  collection of 
complex vector fields $\{X_i\}_{1 \leq i \leq n}$ (with  the same regularity of $J$), which give bases for the $(+i)$-eigenspaces $T^{1,0}_x M$ at all points in which  they are defined. \par
\medskip
\begin{definition}Ê\label{J-holomorphic}ÊA {\it set of $J$-holomorphic coordinates} is a  homeomorphism 
$F = (F^1, \ldots, F^n): \cU \subset M\longrightarrow \cV \subset \bC^n$
 such that \beq \overline{X_i}(F^j) = 0\qquad \text{for any}\ 1 \leq i, j \leq n\eeq
 for any $n$-tuple of  local generators $\{X_i\}_{1 \leq i \leq n}$ of  the distribution $T^{1,0} M\subset T^{\bC} M$ of eigenspaces $T^{1,0}_x M$. 
 \end{definition}
 \medskip
 An atlas $\cA_J$ of smoothly overlapping systems of $J$-holomorphic coordinates makes $(M, \cA_J)$  a complex manifold of dimension $n$. 
 Conversely, if $(M, \cA)$ is a complex manifold, the unique 
 tensor field $J$,   which in all charts of $\cA$ is of the form  
\beq \label{classicalJ} J = i \frac{\partial}{\partial  z^i} \otimes d  z^i -  i \frac{\partial}{\partial \overline{z^i}} \otimes d \overline{z^i}\eeq 
is  a formally integrable almost complex structure, for which the charts of $\cA$ are $J$-holomorphic. \par
\medskip
The celebrated  Newlander-Niremberg Theorem (\cite{NN}) states that a smooth almost complex structure $J$  is  formally integrable if and only if  the set $\cA_J$ of  $J$-holomorphic coordinates  is a complex atlas and  $(M, \cA_J)$ is a complex manifold.
  Therefore on a $2n$-dimensional manifold $M$, there is a  one-to-one correspondence between the  complex manifold structures   and the smooth formally integrable almost complex structures  $J$. Due to this,   any complex manifold can be equivalently identified as the pair  $(M, \cA_{J})$ or   the pair $(M, J)$.  \par
 \par
\medskip
Newlander-Niremberg Theorem follows from the following existence  result, which  spells out the minimal regularity assumptions required for $J$. 
\begin{theo}[\cite{NN,Ma, NW,  We, HT}]\label{theo22} Let $J$ be a formally integrable, almost complex structures on $M$ of class $\cC^{k + \a}$ for some $k +\a > 1$. 
For any $x_o \in M$, there exists a system of $J$-holomorphic coordinates  $F: \cU\longrightarrow \bC^n$ of class 
$\cC^{k+1 + \a}$ on a neighborhood $\cU$ of $x_o$.
\end{theo}
\par
\medskip
Finally we  recall the definition of  the operator $dd^c_J$, which is frequently considered in this paper. For  a given complex structure $J$, 
for any $\cC^2$ function  $f$, one has  
$d^c_J f \= J df  = - df(J (\cdot))$ and
\beq \label{ddc}Êd d^c_{J} f(X, Y) = - X(JY(f)) + Y(JX(f) )+ J[X, Y](f)\ .\eeq
A local computation shows  that
$d d^c_J  = 2 i  \partial \overline \partial$.
\par
\bigskip
\subsection{Domains of circular type and Monge-Amp\`ere exhaustions}\hfill\par
A {\it (bounded) manifold of circular type} 
is a  complex manifold $(M, J)$  endowed  with   a $\cC^0$ exhaustion  $\t: M \longrightarrow [0,1)$ such that
\begin{itemize}
\item[i)]  $\{\t = 0\}$ is a singleton, say  $\{x_o\}$,  Êand  $\t|_{M \setminus\{x_o\}}$ and  the pull back of $\t$ 
 on the blow up $\wt M$ of $M$  at $x_o $ are both of class $\cC^\infty$;
\item[]
\item[ii)] on $M \setminus \{x_o\} = \{\ \t\geq 0\ \}$ 
$$\left\{\begin{array}{l} 2 i \partial \overline\partial \t = d d^c \t > 0\ ,\\
\ \\
2 i \partial \overline \partial \log \t = d d^c \log \t \geq 0\ , \\
\ \\
(d d^c \log \t)^n \equiv 0\ \ (\text{Monge-Amp\`ere Equation)}\ \ ;\end{array}\right.$$
\item[]
\item[iii)]  in some (hence,   {\it any}) system of  complex coordinates $z = (z^i)$ centered at $x_o$, the exhaustion $\t$ has a logarithmic singularity at $x_o$, i.e. 
$$\log \t(z) = \log\|z\| + O(1)\ .$$
\end{itemize}
 We call  $\t$  {\it Monge-Amp\`ere exhaustion}   and $x_o$ {\it center of the exhaustion}. \par
\medskip
Any  Monge-Amp\`ere exhaustion  $\t$ defines an associated distribution   
\beq\label{definitionofZ} Ê \cZ  = \{\ v \in T(M\setminus\{x_o\})\ :\ (\partial \overline \partial \log \t)(v, \cdot ) = 0\ \} \subset T (M \setminus \{x_o\})\ .\eeq
It is  an integrable complex distribution and  its integral leaves are complex curves whose  closures are holomorphic disks passing  through   $x_o$. We  call it  
{\it Monge-Amp\`ere foliation} associated with $\t$.\par
\medskip
We call  {\it domain of circular type} any  pair $(D, \t)$, formed by a relatively compact domain  $D \subset (M, J)$ with  smooth boundary,  
and  an  exhaustion $\t: D \longrightarrow [0, 1)$,  which is smooth  up to the boundary and  such that   $((D,J), \t)$  is a manifold of circular type. \par
\medskip
The   modeling example for such  domains  is  the unit ball $\bB^n \subset \bC^n$, equipped  with the {\it standard Monge-Amp\`ere exhaustion} $\t_o: \bB^n \longrightarrow [0,1)$, $\t_o(y) = \|y\|^2$ and Monge-Amp\`ere foliation   formed  by  the radial  disks  
\beq\label{radialdisk}  \D^{(v)} = \{\ \z v\ ,\  \z \in \D\ \}\ ,\qquad  v \in \bC^n \setminus \{0\}\ .\eeq
The class of domains  of circular type naturally includes all  strictly pseudoconvex  circular domains  and all strictly linearly convex domain of $\bC^n$.   Indeed,  for any strictly pseudoconvex circular domain    $D \subset \bC^n$,  the associated squared Minkowski function $\mu^2_D$ is a  Monge-Amp\`ere exhaustion having $x_o = 0$ as center.   On the other hand, for any 
 strictly linearly convex domain $\O \subset \bC^n$ and  for any choice of point  $x_o \in \O$,  
 the exhaustion 
  \beq \label{lilla} \tau: \O \longrightarrow [0,1)\ ,\qquad \t(y)= (\tanh \delta(y))^{2}\ ,\eeq
  where $\d(y)$ denotes the Kobayashi distance between $y$ and  $x_o$, is  a  Monge-Amp\`ere exhaustion  for $\O$ that has    $x_o$ as center.\par
\bigskip

\subsection{Normal forms of  domains of circular type}\hfill\par
\label{section23bis}
Consider now the domain of circular type $(\bB^n, \Jst, \t_o)$, given by the unit ball $\bB^n$,  with its standard complex structure $\Jst$  and standard exhaustion $\t_o$.  Let also  $\cZ$ the corresponding distribution \eqref{definitionofZ}  on $\bB^n\setminus \{0\}$  and set
\beq \label{distrter}Ê\cH = (\cZ)^\perp = \{\ v \in T_x\bB^n , \ x \neq 0  \ :\ < v,\cZ_x> = 0\ \}\ , \eeq 
where $< \cdot, \cdot>$ denotes  the standard  Euclidean metric.  They  are not defined at $0$
 but they both admit  smooth extensions  at all points of  the blow up $\wt p: \wt \bB^n \longrightarrow \bB^n$ at $0$. 
We call  $\cZ$  {\it (standard) radial distribution} and 
 $\cH$    {\it (standard) normal distribution}.  
The distributions $\cZ$,    $\cH$ are  both $\Jst$-invariant. We denote  $\cZ^{1,0}, \cZ^{0,1} \subset \cZ^\bC$  and $\cH^{1,0}, \cH^{0,1} \subset \cH^\bC$, the complex subdistributions, given  by the $(+i)$- and $(-i)$-eigenspaces of $\Jst$ in $\cZ^\bC$ and $\cH^\bC$, respectively. \par
\medskip
\begin{definition} \label{definition5.1} We call {\it L-complex structure} any formally integrable almost complex structure $\wt J$ on the blow up $\wt \bB^n$  that satisfies the following conditions:
\begin{itemize}
\item[i)] $\wt J$ leaves invariant all spaces of the distributions $\cZ$,  $\cH$; 
\item[ii)] $\wt J|_{\cZ} = \Jst|_{\cZ}$; 
\item[iii)] there exists a smooth homotopy $\wt J_t$ of complex structures on $\wt \bB^n$ that satisfy (i) and (ii), with $\wt J_{t = 0} = \Jst$ and $\wt J_{t = 1}= \wt J$.  
\end{itemize}
Let $\cA_{\wt J}$ be the  atlas on $\wt \bB^n$ of $\wt J$-holomorphic coordinates for an L-complex structure $J$.  The  blow-down   of $(\wt \bB^n, \cA_{\wt J})$ at $0$ is called {\it manifold in normal form associated with $\wt J$}.  A  {\it  manifold (of circular type) in normal form} is any   complex manifold determined in the above fashion. \end{definition}
 \begin{rem} \label{rem24} As   topological space, any manifold $M$ in normal form is  naturally homeomorphic to $\bB^n$. Hence, it can be considered  as a pair $(\bB^n, \cA_J)$,  formed by  $\bB^n$ and an appropriate atlas $\cA_J$ of complex charts with holomorphic overlaps. Here, $J$ denotes for the integrable almost complex structure defined by  \eqref{classicalJ} in the charts of  $\cA_J$. \par
The existence of a blow down complex manifold $(\bB^n, \cA_J)$ for  $(\wt \bB^n, \cA_{\wt J})$ follows directly from   the fact that  $(\wt \bB^n,\cA_{\wt J})$ has  a plurisubharmonic exhaustion, which  is strictly plurisubharmonic outside the singular set $ \wt p^{-1}(0)$  (see e.g., \cite{PS}). However,   the complex charts  of such blow down  are in general  not smoothly overlapping with the standard coordinates of $\bC^n$. This means that  the  tensor field $J$ (which has  smooth components  in the charts of  $\cA_J$) might have non-smooth components when it is expressed in terms of  the standard coordinates of $\bC^n$.   In \S \ref{section31}, we discuss  this  phenomenon in  detail. Here, we just point out that {\it $J$ is surely smooth (in standard coordinates) at the  points of $\bB^n \setminus \{0\}$ and that lower regularity might occur  at  $0$}.
\end{rem}
In  \cite{PS},   it is proved that any manifold in normal form is of circular type and that the following holds:
 \begin{theo} \label{existenceanduniqueness}ÊLet  $(M, \cA_J)$ be  a manifold of circular type  with Monge-Amp\`ere exhaustion $\t$ and center $x_o$.  Then, there exists a biholomorphism 
 $\Phi: (M, \cA_J) \longrightarrow (\bB^n, \cA_{J'})$  with a  manifold in normal form $(\bB^n, \cA_{J'})$ such that 
 \begin{itemize} 
 \item[a)]Ê$\Phi(x_o) = 0$ and $\t = \t_o \circ \Phi$; 
 \item[b)]Ê$\Phi$ maps the leaves of the Monge-Amp\`ere foliation of $M$ into the straight disks through the origin of $\bB^n$. 
 \end{itemize}
 \end{theo}
 A biholomorphism  $\Phi: (M, \cA_J) \longrightarrow (\bB^n, \cA_{J'})$ satisfying   (a) and (b)  is called {\it normalizing map for $(M, \cA_J, \t)$}. \par
 \bigskip
\subsection{Bland and Duchamp invariants of  manifolds of  circular type}\hfill\par
\label{BDinvariants}
  Let   $(\bB^n, \cA_J)$ be a manifold in normal form.  By definition,  the  formally integrable almost complex structure   $J$ on $\bB^n \setminus \{0\}$   differs from the standard complex structure $\Jst$   only by its action on the vectors in  $\cH$.  Hence it is  uniquely determined by the distribution $\cH^{0,1}_J$ of  the  $(-i)$-eigenspaces of $J_x$ in the complex tangent spaces $T^\bC_x \bB^n$, $x \in T_x \bB^n$. In \cite{PS}, it is shown that $\cH^{0,1}_J$  is necessarily  of the form
\beq \label{def1}\left.Ê\cH^{0,1}_{J}\right|_x = \left\{\ v = w+ \left.\phi_J\right|_x(v)\ ,\ v \in \left.\cH^{0,1}\right|_x\ \right\} \subset \ \cH^\bC_x\ ,\qquad x \neq 0\ ,\eeq
for some appropriate  tensor field
$\phi_J \in (\cH^{0,1})^* \otimes \cH^{1,0} $, 
 called  {\it deformation tensor of $J$ with respect $\Jst$}.
 \par
 \medskip
The facts that $J$ is a formally integrable almost complex structure and that   $\t_o = \|\cdot\|^2$  is a  Monge-Amp\`ere exhaustion for $(\bB^n, \cA_J)$ give strong   constraints on the deformation tensor $\phi_J$. In particular, it turns out that $\phi_J$ is  sum of a  series,  uniformly converging  on compacta,   of the form 
\beq \label{sum}Ê\phi_J =  \sum_{k \geq 0}Ê\phi^{(k)}_J\ , \eeq
with each tensor field  $\phi_J^{(k)}$ is  (locally)  of the form
\beq \label{2.9} \phi_J^{(k)}([w], \z w) = \phi^k_J([w]) \z^k\ .\eeq 
Here, we indicate any point  $x \in \bB^n \setminus \{0\} = \wt \bB^n \setminus p^{-1}(0)$  by  the corresponding pair $([w], \z) \in \bC P^{n-1} \times \D$, with $w \in S^{2n-1}$, such that $x = \z \cdot w$. The tensor field $\phi_J$ is  constrained by the equations and  inequality
\begin{itemize}
\item[i)] $dd^c\tau_o(\phi_J(X), Y) + dd^c \tau_o(X, \phi_J(Y)) = 0$ for  any  $X,Y \in \cH^{0,1}$,
\item[ii)] $\bar \partial_b \phi_J + \frac{1}{2}[\phi_J, \phi_J] = 0$, 
 \item[iii)] $ dd^c \tau_o(\phi_J(X), \overline{\phi_J(X)}) < dd^c \tau_o(\bar X, X)$   for any  $0 \neq X \in H^{0,1}$
\end{itemize}
(for the detailed definition of  $\bar \partial_b$, see \cite{PS}). Conversely,   any sequence of tensor fields $\phi_J^{(k)} \in (\cH^{0,1})^* \otimes \cH^{1,0}$, $0 \leq k < \infty$, of the form \eqref{2.9}, such that  \eqref{sum} converges uniformly on compacta of $\wt \bB^n$ and the sum satisfies  (i) - (iii), the corresponding complex distribution \eqref{def1} determines uniquely  an almost complex structure  $J$ on $\wt \bB^n$, 
which  is an $L$-complex structure, hence corresponds to  a manifold in normal form. \par
The terms  $\phi_J^{(k)}$ of the series \eqref{sum} are called {\it Bland and Duchamp invariants} of the manifold  in normal form $(\bB^n, \cA_J)$.\par
\medskip
\begin{lem}  \label{lemmetto}ÊLet $(\bB^n, \cA_J)$ be a manifold  in normal form,  with associated deformation tensor $\phi_J$. If  $X^{1,0}, Y^{0,1}$ are   complex vector fields in  $\cH^{1,0}$,   $\cH^{0,1}$, respectively, then  
\beq \label{lemmettoeq} dd^c_J \t_o(X^{1,0} + \overline{\phi_J(\overline{X^{1,0}})} , Y^{0,1} + \phi_J(Y^{0,1})) =$$
$$ =  d d^c_{\Jst} \t_o(X^{1,0}, Y^{0,1}) +  d d^c_{\Jst} \t_o(\overline{\phi_J(\overline{X^{1,0}})}, \phi_J(Y^{0,1}))\ . \eeq
\end{lem}
\begin{pfns}  Since $[\cH^{1,0}, \cH^{1,0}] \subset \cH^{1,0}$ and $[\cH^{0,1}, \cH^{0,1}] \subset \cH^{0,1}$,
\beq\label{2.10}Ê\left[X^{1,0}\! + \!\overline{\phi\left(\overline{X^{1,0}}\right)}, Y^{0,1}\!+\! \phi\left(Y^{0,1}\right)\right]\!
= 
 \! \left[X^{1,0} , Y^{0,1}\right] 
 + \! \left[ \overline{\phi\left(\overline{X^{1,0}}\right)}, \phi\left(Y^{0,1}\right)\right]\!\!\!\!\!\!\!  \mod\! \cH^{\bC} \eeq
Since $X'{}^{1,0}Ê\=  X^{1,0} + \overline{\phi\left(\overline{X^{1,0}}\right)} \in \cH^{1,0}_J$,   $Y'{}^{0,1}  \= Y^{0,1}+ \phi\left(Y^{0,1}\right) \in \cH^{0,1}_J$ and   the derivatives of $\t_o$ along  vectors in $\cH$ are trivial, we get 
$$Êd d^c_J \t_o(X'{}^{1,0}, Y'{}^{0,1})  \overset{\eqref{ddc}}=  i X'{}^{1,0}\left(Y'{}^{0,1}\left(\t_o\right) \right) - i Y'{}^{1,0}\left(X'{}^{0,1}\left(\t_o\right) \right) - \phantom{aaaaaaaaaaa}$$
$$\phantom{aaaaaaaaaaaaaaaaaaaaaaaaa} - J[X'{}^{1,0}, Y'{}^{0,1}](\t_o) 
=  - J[X'{}^{1,0}, Y'{}^{0,1}](\t_o)  \overset{\eqref{2.10}}=$$
$$ \phantom{aaaaaaaaaaaaaaaa} = - J(\left[X^{1,0} , Y^{0,1}\right])(\t_o)  
 -  J \left(\left[ \overline{\phi\left(\overline{X^{1,0}}\right)}, \phi\left(Y^{0,1}\right)\right] \right) (\t_o)\ .$$
 Now,   using once again that the derivatives of $\t_o$ along  vectors in $\cH$ are trivial and since    $J|_{\cZ} = \Jst|_{\cZ}$, 
 $$Êd d^c_J \t_o(X'{}^{1,0}, Y'{}^{0,1}) =  - \Jst(\left[X^{1,0} , Y^{0,1}\right])(\t_o)  
 -  \Jst \left(\left[ \overline{\phi\left(\overline{X^{1,0}}\right)}, \phi\left(Y^{0,1}\right)\right] \right) (\t_o) = $$
 $$ =  d d^c_{\Jst} \t_o(X^{1,0}, Y^{0,1}) +  d d^c_{\Jst} \t_o(\overline{\phi_J(\overline{X^{1,0}})}, \phi_J(Y^{0,1}))\ .\eqno{\square}$$
\end{pfns}
\bigskip
 \subsection{Kobayashi indicatrices and normalizing  maps}\hfill\par
 \label{section25}
 Let $M$ be a manifold of circular type with Monge-Amp\`ere exhaustion $\t$ and center $x_o$.  We recall  that  the value $\k(v)$ of the infinitesimal Kobayashi metric $\k$ at  any vector $v \in T_{x_o} M$ always exists and can be computed by the expression  
$Ê\k(v) = \lim_{t \to 0}Ê\left.\frac{d}{dt} \sqrt{\t(\g_t)}\right|_{t = t_o}$, 
where   $\g_t$ is any  smooth curve  such that $\g_0 =x_o$ and $\dot \g_0 = v$ (\cite{Pt3}). The  {\it Kobayashi indicatrix of $M$ at $x_o$} is the    domain $I_{x_o} $ in the tangent space $ T_{x_o} M $,  defined by 
$$I_{x_o} = \{\ v \in T_{x_o}M\ :\ \k(v) < 1\ \}\ .$$
Given a basis $\cB = (e_i)$ for $T_{x_o} M$,   we  denote by $\ell_{\cB}$ the $\bC$-linear isomorphism 
$$\ell_{\cB}: \bC^n \longrightarrow T_{x_o} M\qquad \text{such that}\qquad  \ell_\cB(v^i e^o_i) \= v^i e_i\ .$$
The  circular domain of $\bC^n$ corresponding to $I_{x_o}$
\beq \label{idenfiedindicatrix}ÊI^{(\cB)} \= \ell_\cB^{-1}(I_{x_o}) =  \{\ v \in \bC^n\ : \ \k(\ell_\cB(v)) < 1 \ \} \ .\eeq
is called   {\it realization of  $I_{x_o}$ determined by the basis  $\cB$}.
\par
Let $p: \wt I^{(\cB)} \longrightarrow I^{(\cB)}$
be  the blow up of the realization  $I^{(\cB)}$ at  $0 \in \bC^n$. We recall that $\wt I^{(\cB)}$ is naturally identifiable with the  set of  pairs $([v], t v)$,  formed by 
elements $[v]Ê\in \bC P^{n-1}$, with $v \in \partial I^{(\cB)}$,  and  the points $t v \in I^{(\cB)}$, $t \in [0,1)$. 
 Let  also     $\pi: \wt M \to M_o$ be  the blow up of $M$ at $0$ and $x_o$. \par
\medskip
 The {\it circular representation}  of the manifold of circular type  $(M,\cA_J)$ is a   diffeomorphism 
 $$ \Psi: \wt I^{(\cB)}Ê\longrightarrow \wt M \ , $$
  which is canonically  determined by the Monge-Amp\`ere exhaustion  $\t$ and satisfies a number of crucial properties. For  the definition of $\Psi$ and its main properties, we refer to e.g. \cite{Pt, Pt1,  BDK, PS}. Here, we just remind that:   
\begin{itemize}
\item[a)] for any  $([v],  t v) \in \wt I^{(\cB)}$,  $v \in \partial I^{(\cB)}$,   the  map 
$$f^{v}: \D \longrightarrow M\ ,\qquad f^{v}(\z) \= (\pi \circ  \Psi)([v], \z v) $$
is the unique stationary disk of $M$  satisfying the conditions
$f^{v}(0) = x_o$ and  $f^{v}_*\left(\left.\frac{\partial}{\partial x}\right|_{0}\right) = \l v$ for some $\l > 0$; 
moreover, the  corresponding  lifted disk 
$\wt f^{v}: \D \longrightarrow \wt M$, defined by  $\wt f^{v}(\z) \= \Psi([v], \z v)$
is proper, holomorphic and injective; 
\item[b)]Êthe restriction 
 $\Psi|_{\pi^{-1}(0)} : \bC P^{n-1}Ê\longrightarrow \pi^{-1}(x_o) = P(T_{x_o}M) \subset \wt M$ coincides with the projective map
 $\wt  \ell_{\cB}$, determined by $\ell_B$, i.e., 
$\Psi|_{\pi^{-1}(0)}([v]) = \wt \ell_\cB([v]) \= [\ell_\cB(v)]$. 
 \end{itemize}
\par
 \smallskip
Moreover by  \cite{PS}, Lemma 3.5, there exists a diffeomorphism   $\psi: \wt \bC^n \longrightarrow \wt \bC^n$   from the blow up of $\bC^n$ at the origin into itself satisfying: 
 \begin{itemize}
 \item[1)]Êit is a bundle map w.r.t. the   $\bC^*$-bundle structure $\wt \pi: \wt \bC^n \longrightarrow  \bC P^{n-1}$;  
 \item[2)] it preserves the distribution $\cH$ of $\wt \bC^n$, formed by the  vectors that are normal to  the  straight  lines of $\bC^n$  through the origin; 
 \item[3)] it  induces a diffeomorphism between   $\wt I^{(\cB)}Ê$  and   $\wt \bB^n$; 
 \item[4)]Êthere exists a smooth homotopy between $\psi$ and $\Id_{\wt \bC^n}$,  for which   all intermediate maps $\psi_t$ satisfy (1) and (2) and map $\wt I^{(\cB)}$ onto a blow up at the origin of a circular domain. 
 \end{itemize}
\par
 \smallskip
 It turns out that the diffeomorphism 
 $\wt \Phi \= \psi \circ \Psi^{-1}: \wt M^n \longrightarrow \wt \bB^n$ 
 maps the integrable almost complex structure $\wt J$ of $\wt M^n$ into the almost complex structure $\wt J' = \wt \Phi_*(\wt J)$ on $\wt \bB^n$, which is 
 an $L$-complex structure,  and  $(\wt \bB^n, \cA_{\wt J'})$ is the blow up of a manifold in normal form $(\bB^n, \cA_{J'})$.  Hence, the  pushed-down diffeomorphism  $\Phi: (M, \cA_J) \to (\bB^n, \cA_{J'})$ is a normalizing map for  $(M, \cA_J, \t)$.\par 
\smallskip
As a  consequence of the above construction of a normalizing map,   we get:
\begin{lem} \label{lemma27} If there exists a  basis $\cB$ for $T_{x_o}ÊM$ such that   $I^{(\cB)} =  \bB^n$,  there also  exists a normalizing map $\Phi: (M, \cA_J) \longrightarrow (\bB^n, \cA_{J'})$ onto a manifold in normal form whose    deformation tensor  $\phi_{J'}$ has  $\phi^{(0)}_{J'} = 0$. 
\end{lem}
\begin{pf} If  $I^{(\cB)} =  \bB^n$,  as a diffeomorphism $\psi$ satisfying (1)-(4), we may take  the identity  map $\psi =  \Id_{\wt \bC^n}$.  Hence,  the inverse $\wt \Phi = \Psi^{-1}: \wt M \longrightarrow \wt I^{(\cB)} = \bB^n$, of the circular representation   of $M$   is a normalizing map. By  (b) and the fact that $\ell_\cB$ is a $\bC$-linear isomorphism,  we get that 
$$\wt J|_{\bC P^{n-1}} = \Psi^{-1}_*(J|_{P(T_{x_o} M)}) =  \wt \ell_{\cB}^{-1}{}_*(J|_{P(T_{x_o} M)}) = \Jst$$
where  $\Jst$ denotes  the standard complex structure of  $\bC P^{n-1}$. This implies  that the deformation tensor $\phi_{\wt J}$ of $\wt J$ vanishes  on $\bC P^{n-1}$ and    $\phi^{(0)} = 0$. 
\end{pf}
\bigskip
\section{The almost complex structures of   manifolds in normal form }
\subsection{The tangent bundle  of  a manifold in normal form}\hfill\par
\label{section31} 
\setcounter{equation}{0}
 Let  $(\bB^n, \cA_J)$ be a  manifold in normal form and denote by $\cA^{\bR}_J$  the complete $\cC^\infty$ atlas determined by the real coordinates, given by real and imaginary parts of the complex coordinates in $\cA_J$. \par
By Remark \ref{rem24}, $\cA^{\bR}_J$  is in general different from the standard atlas  $\cA^{\bR}_o$ of $\bB^n$, induced by the   standard  real manifold structure   of   $\bR^{2n} (= \bC^n)$.\par
\smallskip
A priori,    the   ``tangent vectors of $\bB^n$''   are  different from the usual  vectors of $\bB^n$  if  we  consider the non-standard  atlas $\cA^{\bR}_J$ on  $\bB^n$ in place of the standard one.   Indeed,  the tangent vectors of $(\bB^n, \cA^{\bR}_o)$  are   equivalence classes of curves, which are  of class   $\cC^1$ when they are expressed in some coordinates of the standard atlas. On the contrary, the tangent vectors of $(\bB^n, \cA^{\bR}_J)$ are equivalence classes of curves, which are  of class $\cC^1$ in the charts of the atlas $\cA^{\bR}_J$ and  might not be of class $\cC^1$ for the charts in  $\cA^{\bR}_o$. \par
\medskip
Nonetheless, as we shall see below, {\it there  is a  natural identification between the   tangent vectors of $(\bB^n, \cA^{\bR}_o)$ and  $(\bB^n, \cA^{\bR}_J)$}.  Using this fact, in the next sections, we  shall not make any distinction between the tangent spaces of such two differentiable manifolds.  
\par
\bigskip
Consider  the blow-up $(\wt \bB^n,  \cA_{\wt J})$ of $(\bB^n, \cA_J)$ at the origin. By construction,  
$\wt J$  is a formally integrable,  almost complex structure, which is smooth in the standard atlas of real coordinates $\cA_o^{\bR}:= \cA^{\bR}_{\wt J_o}$ of the blow up $(\wt \bB^n, \wt \Jst)$ of $(\bB^n, \Jst)$.  Hence, by  Theorem \ref{theo22}, every complex chart   in $\cA_{\wt J}$ (hence, every real chart in  $\cA^{\bR}_{\wt J}$)   smoothly overlaps with all   charts  of  $\cA_o^{\bR}$. 
This means that the   atlas $\cA^{\bR}_{\wt J}$  coincides with the  standard   atlas  $\cA^{\bR}_o$ of $\wt \bB^n$. 
Note also that  the blow-down maps $\wt \pi: (\wt \bB^n, \cA_{\wt \Jst}) \longrightarrow (\bB^n, \cA_o)$ and $\wt \pi': (\wt \bB^n,  \cA_{\wt J}) \longrightarrow (\bB^n, \cA_J)$,  as maps between the underlying topological manifolds $\wt \bB^n$ and $\bB^n$, do coincide: $\wt \pi = \wt \pi'$. 
Furthermore the restriction $\wt \pi|_{\wt \bB^n \setminus \wt \pi^{-1}(0)}=
\wt \pi'|_{\wt \bB^n \setminus \wt \pi'{}^{-1}(0)}$ 
 is a diffeomorphism (indeed, it is a biholomorphism w.r.t. the complex structures involved). \par
 \medskip
 These observations show that   if a  chart of  the atlas $\cA_J$   is defined on some open subset of  $\bB^n \setminus \{0\}$,   it smoothly overlaps with all charts of $\cA_o$. This implies that that   curves  in $\bB^n \setminus \{0\}$    are of class $\cC^1$  in the charts of  $\cA_J$ if and only if   are of class $\cC^1$ in the charts of  the atlas $\cA_o$ and vice versa. 
 Hence  {\it  any tangent space 
 $T_x \bB^n$, $x \neq 0$, of the real  manifold  $(\bB^n, \cA_J)$ is actually identical  with  the   tangent space  at $x$ of the standard ball  $(\bB^n, \cA_o)$}. \par
 \medskip
Let us now focus on   the tangent spaces  at $x = 0$  of   $(\bB^n, \cA_J)$ and $(\bB^n, \cA_o)$. 
A priori,  the codomain of the push-forward   $\pi_*: T \wt \bB^n \longrightarrow T\bB^n$
changes if     one  considers $T\bB^n$ as tangent bundle of $(\bB^n, \cA_J)$ or as  tangent bundle of   $(\bB^n, \cA_o)$. However, in both cases,  the tangent space $\bB^n$ at the origin is equal to the image $ \pi_*(T\wt \bB^n|_{\pi^{-1}(0) })$. This brings to a natural   identification between the two tangent spaces, which we now explain.
\par
Firstly we observe  that if  $v$, $v' \in T\wt \bB^n|_{\pi^{-1}(0) }$ are such that 
$ \pi_*(v) = \pi_*(v') $ when the target manifold is $(\bB^n, \cA_o)$ 
 then one also has
$ \pi_*(v) = \pi_*(v')$ when the target manifold is $(\bB^n, \cA_J)$. 
This can be checked as follows. When the target manifold is $(\bB^n, \cA_o)$, 	we have $ \pi_*(v) = \pi_*(v')$  exactly if  there is a pair  of smooth curves $\wt \g, \wt \g' : (-\varepsilon, + \varepsilon) \to \wt \bB^n$, 
 with
 \begin{itemize}
 \item[i)] $y = \wt \g(0)$ and $y' = \wt \g'(0)$ are both in $ \pi^{-1}(0)$ and    tangent to $v$ and $v'$, respectively; 
 \item[ii)]  the  curves $\g = \pi \circ \wt \g$ and $\g' = \pi \circ \wt \g'$ in the manifold $(\bB^n, \cA_o)$ are tangent one to the other at $0$. 
 \end{itemize}
 Since the differentiable structures induced on $\bB^n \setminus \{0\}$ of $(\bB^n, \cA_o)$ and $(\bB^n, \cA_J)$ are equal, by a direct argument based on the continuity of tangent vectors,
 the curves  $\g = \pi \circ \wt \g$ and $\g' = \pi \circ \wt \g'$ are tangent   at $t=0$ also when  considered as  curves of  $(\bB^n, \cA_J)$. This proves the claim.\par
\smallskip
By this observation,  {\it we may  identify  every vector $v$ of  the tangent space  $T_0 \bB^n$  of $(\bB^n, \cA_J)$ with a unique corresponding  vector  $\check v$ of the tangent space  $T_0 \bB^n$  of $(\bB^n, \cA_J)$, provided that   there exists some $\wt v \in T\wt \bB^n|_{ \pi^{-1}(0)}$,  which is projected  onto $v$ and $\check v$, respectively, by the corresponding push-forwards}. \par
 \begin{rem} \label{remark31}Ê
With the above identifications between  tangent bundles,  the tensor field $J$ of the manifold $(\bB^n, \cA_J)$ can be considered    as a (possibly non-smooth) tensor field of the standard  ball  $(\bB^n, \cA_o)$. 
 We claim such tensor field satisfies 
 \beq \label{3.3}ÊJ|_0 = \Jst|_0\ .\eeq
  Indeed, $J|_0$ is uniquely determined  by its  $(+i)$-eigenspace in $T_0^\bC \bB^n\ $.
   On the other hand,  such eigenspace is spanned by the tangent vectors of the stationary disks of $(\bB^n, \cA_J)$ passing  through $0$.
 This implies  \eqref{3.3} as 
such  disks  of  $(\bB^n, \cA_J)$ coincide with the stationary disks of $(\bB^n,\cA_o)$ through  $0$.
 \end{rem}

\medskip
\subsection{Almost complex structures of manifolds in normal form} \hfill\par
Let $(\bB^n, \cA_J)$ be a manifold  in normal form, with formally integrable almost complex structure $J$. Using Remark \ref{remark31}, we may consider $J$ as a tensor field on the  standard unit ball  $(\bB^n, \cA_o)$. 
 \par
 \smallskip
As mentioned before, the components of    $J$  in the  charts of $\cA_J$ are smooth, but, a priori, they are not smooth in the standard coordinates. 
Actually,   since the charts of $\cA^{\bR}_J$ on open subsets of $\bB^n \setminus \{0\}$,  smoothly overlap with the standard coordinates, the components of $J$ at  points of $\bB^n \setminus \{0\}$ are necessarily smooth also in standard coordinates.  On the other end, the  components  of $J$ at $0$ in standard coordinates  might  not  even be  continuous. \par
\medskip
 From now on, if no indication of the atlas is given,  {\it  $\bB^n$ is considered as endowed with  charts of the standard atlas $\cA_o$ and  $J$ is considered as a formally integrable,  almost complex structure on $(\bB^n, \cA_o)$,  which is  of class $\cC^\infty$ at all points of $\bB^n \setminus \{0\}$}. \par
 \medskip
 By construction,  the Nijenhuis tensor $N_J$ of $J$ vanishes identically at all points of $\bB^n \setminus \{0\}$, i.e.  $J|_{\bB^n \setminus \{0\}}$ is   formally integrable. 
 Moreover, every  system of $J$-holomorphic  coordinates of the (non-standard) manifold $(\bB^n, \cA_J)$
 $$\xi_J = (z^1_J, \ldots, z^n_J): \cU \longrightarrow \bC^n$$
 is a system of $J$-holomorphic coordinates for the almost complex structure $J$ of the (standard) manifold $\bB^n$, as defined in  Definition \ref{J-holomorphic}. \par
 \medskip
Hence {\it for any point  $x_o \in \bB^n$ there exists a system of $J$-holomorphic coordinates $\xi_J = (z^i_J)$ defined on a neighborhood of $x_o$}.  This is {\it not}  a consequence of the general Theorem \ref{theo22} as $J$ might be not even $\cC^0$ at $0$. It simply follows from the fact that  $J$ corresponds to a manifold in normal form, namely  to $\bB^n$  endowed with a (non-standard) complex manifold structure. \par 
 \medskip
 \begin{lem} Let $(\bB^n, J)$ be a manifold  in normal form. Then
 \begin{itemize}
 \item[a)] if there exists a 
 system of $J$-holomorphic coordinates $\xi_J = (z^j_J)$  in a  neighborhood $\cU$ of $0$,
which is of class $\cC^{1 + r + \a}$ at $0 \in \bB^n$, for an integer $r\geq 0$ and $\a \in [0,1)$, then 
  $J$ is  of class $\cC^{r+\a}$   at the origin.  
 \item[b)]ÊIf  $J$ is of class $\cC^{r +\a}$ at $0$,  for an integer $r \geq 1$ and $\a \in (0,1)$, then there exists a 
 system $\xi_J = (z^j_J)$    of $J$-holomorphic coordinates  of class $\cC^{1 + r + \a}$ in a  neighborhood $\cU$ of $0$. \par
\end{itemize}
 \end{lem}
 \begin{pf} (a) Assume that  $\xi_J = (z^j_J)$ is of class $\cC^{1+r + \a}$ at $0$. This implies that  the entries of the Jacobians
 $$\left.\left(\begin{array}{cc}Ê \frac{\partial z^i_J}{\partial z^\ell} \ & \ \frac{\partial z^i_J}{\partial \overline z^\ell} \\
 \ &\ \\
 \frac{\partial \overline z^i_J}{\partial z^\ell} &  \frac{\partial \overline z^i_J}{\partial \overline z^\ell}\end{array}\right)\right|_{(z^m)}\ ,\qquad 
 \left.\left(\begin{array}{cc}Ê \frac{\partial z^i}{\partial z^\ell_J} \ & \ \frac{\partial z^i}{\partial \overline z^\ell_J} \\
 \ &\ \\
 \frac{\partial \overline z^i}{\partial z^\ell_J} &  \frac{\partial \overline z^i}{\partial \overline z^\ell_J}\end{array}\right)\right|_{(z^\ell_J(z^m))}$$
 are $\cC^\infty$ at all points $(z^m) \neq 0$ and of class $\cC^{r+\a}$ at $0$. Using the transformation rules of (real) tensors of type $(1,1)$ and the fact that $J$ has constant components in the $J$-holomorphic  coordinates  $\xi_J = (z^i_J)$, the components of $J$ in  standard coordinates  have the same regularity of the entries of the above Jacobians. In particular, they are of class $\cC^{r+\a}$ at $0$.  \par
 (b) follows directly from Theorem \ref{theo22}
  \end{pf}
\bigskip
\section{The vanishing of the first Bland and Duchamp invariants\\ implies
 regularity  of the  Monge-Amp\`ere exhaustion}
\setcounter{equation}{0}

\subsection{Complex polar coordinates}\hfill\par
As usual,   we indicate the points of the blow up $\pi: \wt \bB^n \longrightarrow \bB^n$  by  pairs $([v],  r v)$,   where $v \in  S^{2n-1}Ê\subset T_0 \bB^n$, $r \in [0,1)$ and   $[v]$ is the equivalence class of $v$ in  $P(T_0 \bB^n \setminus \{0\})  \simeq  P(\bC^n  \setminus \{0\}) = \bC P^{n-1}$.   Recall that  $\wt \bB^n$ can be also considered as a  bundle over $\bC P^{n-1}$,  with projection $\wt \pi:  \wt \bB^n \longrightarrow \bC P^{n-1}$ defined by  $\wt \pi([v], r v) \= [v]$, and standard fiber $\D$.\par
\smallskip
For a fixed $1 \leq i_o \leq n$, we denote by  $ (w^1, \ldots, w^{n-1})$  the usual
  affine  coordinates on $\bC P^{n-1} \setminus \{ [z]\ : \ z^{i_o} = 0\}$ 
$$ ([z^1:  \ldots: z^{n-1}: z^n]) \longrightarrow \left( w^1 \= \frac{z^1}{z^{i_o}}, \ldots \underset{i_o}{\wh{\phantom{i_o}}}Ê\ldots,w^{n-1} \= \frac{z^{n}}{z^{i_o}}\right)$$
and we  call {\it   complex polar coordinates for $\wt \bB^n$} the  associated   coordinates
$\wt \xi_{i_o} = (w^1, \ldots, w^{n-1}, w^n)$ on $\wt \cU_{i_o} =  \pi^{-1}(\{z^{i_o}  \neq 0\})  \subset \wt \bB^n$
defined by 
\beq\label{polar} ([v], r v) \overset{\wt \xi_{i_o} } \longrightarrow  (w^1([v]), \ldots, w^{n-1}([v]),  w^n \= r e^{{\rm i}\arg(z^{i_o})})\ ,\eeq 
for any $ v = (z^1, \ldots, z^n) \in S^{2n-1}$
(here,   $ \arg(\z)$ denotes  the  real number in  $ [0, 2\pi)$ such  $\z = |\z|e^{{\rm i}  \arg(\z)}$). We finally call {\it  complex polar coordinates for $\bB^n$}  the  coordinates $\xi_{i_o}$ on the set $\bB^n \setminus \{ z^{i_o} = 0\}$, defined by 
$$\xi_{i_o}(z) \=   ( w^1(\pi^{-1}(z)), \ldots, w^{n-1}(\pi^{-1}(z)), w^n(\pi^{-1}(z)) )\ .$$
For our purposes, it is  useful  to have explicit  information on  the Jacobians
\beq \label{jacobians} \left.\left(\begin{array}{cc} \frac{\partial w^i}{\partial z^j} & \frac{\partial w^i}{\partial \bar z^j} \\
\ &\ \\
\frac{\partial \bar w^i}{\partial z^j} & \frac{\partial \bar w^i}{\partial \bar z^j}
\end{array} \right)\right|_{(z^i)}\ ,\qquad \left.\left(\begin{array}{cc} \frac{\partial z^i}{\partial w^j}  &  \frac{\partial z^i}{\partial \bar w^j} \\
\ &\ \\
\frac{\partial \bar z^i}{\partial w^j}  &  \frac{\partial \bar z^i}{\partial \bar w^j} 
 \end{array} \right)\right|_{w^m(z^j)} 
\eeq
of  changes of coordinates between  standard  coordinates $ (z^i)$   of $\bB^n$ and    
 complex polar coordinates $ (w^i) = (w^a,  \z)$, $1 \leq a\leq n-1$,   on 
 $ \bB^n \setminus \{ z^{i_o} = 0\}$ for some $1 \leq i_o \leq n$.   We show how to compute them   for the case $i_o = n$:  the other cases can be handled in the same way. \par
 \smallskip
 Let us recall that   if a point has complex polar coordinates  $(w^i) = (w^a, \z) $,  the corresponding standard  coordinates are 
$$z^a = \frac{\z w^a}{\sqrt{1 + \sum_{b = 1}^{n-1} |w^b|^2}} = \mu_w  w^a \z \ ,\qquad z^n =  \frac{\z}{\sqrt{1 + \sum_{b = 1}^{n-1} |w^b|^2}} = \mu_w \z\ ,$$
where  $\mu_w := \frac{1}{\sqrt{1 + \sum_{b = 1}^{n-1} |w^b|^2}}$. 
It follows that the second Jacobian in \eqref{jacobians}  is  
$$  |\z|
\left( \begin{array}{c|c||c|c} \smallmatrix   (\d_a^b -  \frac{1}{2}Ê\mu_w^2  w^b \overline w^c \d_{ca})\mu_w e^{{\rm i} \q}\endsmallmatrix  & \smallmatrix \endsmallmatrix \smallmatrix w^b \mu_w \frac{1}{|\z|}  
\endsmallmatrix & \smallmatrix - \frac{1}{2}\  w^b w^c \d_{ca}  \mu_w^3 e^{{\rm i} \q}\endsmallmatrix & 0\\
& & & \\
\hline &  &  & \\
\smallmatrix - \frac{1}{2}  \overline{w^c} \d_{ca} \mu^3_w e^{{\rm i}\q} 
\endsmallmatrix &\smallmatrix  \mu_w \frac{1}{|\z|}\endsmallmatrix & \smallmatrix - \frac{1}{2}   w^c \d_{ca} \mu_w^3 e^{{\rm i} \q}
\endsmallmatrix & 0 \\
& & & \\
\hline \hline & & & \\
\smallmatrix - \frac{1}{2}  \overline{w^b} \overline{w^c} \d_{ca}\mu^3_w  e^{-i \q}  \endsmallmatrix  & 0 &\smallmatrix    (\d_a^b -  \frac{1}{2}Ê\mu_w^2 \overline{w^a}  w^c \d_c^b) \mu_w e^{-i \q} \endsmallmatrix & \smallmatrix \overline w^b \mu_w \frac{1}{|\z|}\endsmallmatrix \\
& & & \\
\hline
& & & \\
 \smallmatrix - \frac{1}{2}\overline{w^c} \d_{ca} \mu^3_w e^{-i \q} 
\endsmallmatrix  & 0  
&\smallmatrix - \frac{1}{2} w^c \d_{ac} \mu^3_w e^{-i \q}  
\endsmallmatrix &
\smallmatrix  \mu_w \frac{1}{|\z|}\endsmallmatrix 
\end{array} \right)$$
where  $e^{{\rm i} \q} \= \frac{\z}{|\z|}$. The  first  Jacobian in \eqref{jacobians} is clearly the inverse  of this  matrix.   Computing it explicitly  and recalling that $|\z| = \|z\|$ one can  see that 
\beq \label{secondJ} \left.\left(\begin{array}{cc} \frac{\partial w^i}{\partial z^j} & \frac{\partial w^i}{\partial \bar z^j} \\
\ &\ \\
\frac{\partial \bar w^i}{\partial z^j} & \frac{\partial \bar w^i}{\partial \bar z^j}
\end{array} \right)\right|_{(z^i)}= \frac{1}{\|z\|} A(z, \bar z)\qquad \text{for all}\  z \in \bB^n \setminus \{ z^{i_o} = 0\}\ ,\eeq 
where  $A(z, \bar z)$ is a matrix with smooth and uniformly bounded entries. \par
\bigskip
\subsection{The vanishing of   first   $\phi^{(j)}_J$'s and the regularity of $J$}\hfill\par
Let   $(\bB^n, \cA_J)$ be a manifold   in normal form  with associated deformation tensor $\phi_J =  \sum_k \phi^{(k)}_J$.
Denote by  $\Pi$  the  tensor field of type $(1,1)$ on $\bB^n \setminus \{0\}$, which  gives  the natural projections 
$\Pi_x: T_x \bB^n =\cZ_x + \cH_x  \longrightarrow \cH_x$
 and let   $\psi_J$ be the real tensor field of type $(1,1)$ on $\bB^n \setminus \{0\}$, given by 
\beq \label{psidefinition} \psi_J(X) \=  \phi_J( \Pi(X))\ ,\ \psi_J(\overline{X}) \= \overline{\psi(X)} \ ,\quad \text{for any}\  X\ \text{in}\ T^{0,1} (\bB^n\setminus \{0\})\ .\eeq
  By construction,  
$J|_{\bB^n \setminus \{0\}} = \Jst|_{\bB^n \setminus \{0\}} + \psi_J|_{\bB^n \setminus \{0\}}$.
 On the other hand, by Remark \ref{remark31}, 
$J|_0 = \Jst|_0$. Hence, if we set $\psi_J|_0 \= 0$,  we  have the equality 
$$J = \Jst + \psi_J$$
  at all points of $\bB^n$ and  the regularity of 
 $J$ at $0$  is  equal to  the regularity  of $\psi_J$ at  that point.
Being   the  2-form $d d^c \t_o$   
 non-singular and smooth at all points,   such  regularity   is in turn equal to 
  the regularity at $0$ of the $(0,2)$-tensor field
 \beq \label{tildepsi} \wt \psi_J \= d d^c \t_o(\cdot, \psi_J(\cdot))\ .\eeq
\begin{theo}\label{prop41} Let $k \geq 3$.  If  all  Bland and Duchamp invariants  $\phi^{(j)}$,   $0 \leq j \leq k-1$,     are identically vanishing, then
\begin{itemize}
\item[a)] the almost complex structure $J$ is of class $\cC^{\left[\frac{k}{2}\right] - 1}$  at the origin;  
\item[b)] if  $k \geq 6$,   the standard exhaustion $\t_o: \bB^n \longrightarrow [0,1)$   is  of class $\cC^{\left[\frac{k}{2}\right]}$ in any set of  coordinates of the  atlas $\cA_J$ .
\end{itemize}
\end{theo}
\begin{pf} Consider the components $\wt \psi_{ij}$ of the tensor field \eqref{tildepsi} in a system of complex polar coordinates $ (w^a, \z)$ on an open subset $\cU$ of $\bB^n \setminus \{0\}$, with $0 \in \partial \cU$.  By the vanishing of  $\phi^{(j)}$,   $0 \leq j \leq k-1$, and  definition of $\wt \psi$,  they are  of the form 
${\wt \psi_{ij}}(w^a, \z, \bar w^a, \bar \z)Ê= |\z|^{k} g_{ij}(w^a,\z,  \bar w^a, \bar \z)$ 
for some smooth, uniformly bounded  function $g_{ij}(w^a, \z, \bar w^a, \bar \z)$.
 From \eqref{secondJ} and the transformation rules of the coordinate components of  $(0,2)$ tensors,  it follows that the components of $\wt \psi_J$ in  standard 
coordinates  are  functions of the form $ \|z\|^{k-2} h_{ij}(z, \bar z)$,  
for some  functions $h_{ij}(z, \bar z)$ on $\cU$  that are  uniformly  bounded and with uniformly bounded derivatives. 
Since $0\in \partial\cU$ for any $\cU$ where complex polar coordinates are defined,  if $k \geq 3$, the functions $\wt \psi_{ij}$  are of class $\cC^{\left[\frac{k}{2}\right] - 1}$ at the origin and,  consequently, the same holds   for  $\wt \psi$ and $J$. Furthermore, if    $k \geq 6$,    Theorem \ref{theo22} applies and     any  system of $J$-holomorphic coordinates $(z^i_J) \in \cA_J$ is of class $\cC^{\left[\frac{k}{2}\right] }$ with respect to the standard coordinates. Being smooth in  standard coordinates,  we get  that $\t_o$  is   of class   $\cC^{\left[\frac{k}{2}\right] }$   in  any  chart   of   $\cA_J$.  
\end{pf}
\bigskip

\section{Regularity of the Monge-Amp\`ere exhaustion implies \\
 that the first  Bland and Duchamp invariants vanish}
\setcounter{equation}{0}

\subsection{$\cC^2$-regularity of the  exhaustion and   vanishing of    $\phi^{(0)}_J$}\label{52}\hfill\par

\begin{lem} \label{lemma1.1} Let $(M, J)$ be a manifold of circular type with Monge-Amp\`ere exhaustion $\t$ and center $x_o$. If  $\t: M \longrightarrow [0, 1)$ has 
second order   derivatives  at $x_o$,  there exists a basis $\cB = (e_i)$ of  $T_{x_o} M $ such that (the   realization  Êof)  the indicatrix  at $x_o$   is 
$I^{(\cB)}Ê= \bB^n$.
\end{lem}
\begin{pf} By  \cite{Pt3},  for any $v \in   T_{x_o} M \setminus \{0\}$ there is a unique  Kobayashi extremal disk $f^{v}: \D \longrightarrow (M,J)$  with  $f(0) = x_o$ and $\Re\left(f^{v}\left(\left.\frac{\partial}{\partial \z}\right|_0) \right)\right) \in \bR v$.  This disk is  such that  
$\t(f^{v}(\z)) = |\z|^2$ for any $\z \in \D $ and 
\beq\label{ulla} \k(v) = 1\qquad \text{if and only if} \qquad \Re\left(f^{v}_*\left(\left.\frac{\partial}{\partial \z}\right|_0\right) \right) = v\ . \eeq
It $\t$ has second order derivatives  at $x_o$,  the $2$-form $d d^c_J\t$ is well defined at $x_o$ and,  for any $v \neq 0$, 
\beq\label{ullalla} 1 = \left.\partial \overline \partial (|\z|^2) = - 2i d d^c_{\Jst}( |\z|^2) \right|_{\z = 0}=  - \left.2i d d^c_J \t\right|_{x_o}\!\!\left(f^{v}_*\left(\!\frac{\partial}{\partial \z}\right), f^{v}_*\left(\frac{\partial}{\partial \bar \z}\right)\!\right) = $$
$$ = \frac{1}{\k(v)^2} \left(- 2i \left.d d^c_J \t\right|_{x_o}(v,v)\right) \ .\eeq
 This shows that the  $J$-Hermitian tensor  $- 2 i d d^c_J \t|_{x_o}$ is positively defined. Consider a  basis  $\cB = (e_1, \ldots, e_n)$  for $T_{x_o} M$,  which is unitary  with respect to $- 2 i d d^c_J \t|_{x_o}$. Let  $\ell_\cB: \bC^n \longrightarrow T_{x_o} M$ be 
 the    $\bC$-linear map  with  $\ell_\cB(e^i_o) = e^i$. By construction, 
 $\ell_\cB{}_*(- 2 i d d^c \t|_{x_o}) = < \cdot , \cdot> $
 and $I^{(\cB)} = \ell^{-1}_\cB(I_{x_o}) = \bB^n$.
\end{pf}
As a direct corollary of  Lemma \ref{lemma1.1} and Lemma \ref{lemma27}  we have the following:
\begin{theo} \label{theorem53} Let  $(M, J)$ be a manifold of circular type with Monge-Amp\`ere exhaustion $\t$ and center $x_o$. If $\t$ is of class $\cC^2$ at $x_o$,   the Kobayashi indicatrix  at $x_o$ is linearly equivalent to $\bB^n$ and  $(M, J)$ is biholomorphic, via a normalizing map,  to a manifold in normal form $ (\bB^n, J')$ with   $\phi^{(0)}_{J'} = 0$. 
\end{theo}

From the results of \S \ref{sei},  it will be clear that there is an infinite dimensional family of such manifolds, none of them biholomorphic to  $(\bB^n, \Jst)$.
\par
\bigskip
\subsection{$\cC^{2k}$-regularity  of  the exhaustion and   vanishing of    $\phi^{(j)}_J$, $j\leq  k-1$}\hfill\par
Let   $(M,   J)$ be  a manifold of circular type   with a Monge-Amp\`ere exhaustion $\t: M \longrightarrow [0,1)$ and  center $x_o$. We observe  that
the  2-form $\o = d d^c_J \t$  is a K\"ahler form  for the complex manifold $(M \setminus \{x_o\}, J)$. Furthermore,  
\begin{prop}[\cite{Bu}, Prop.1.1] \label{prop53}The integral  leaves of the distribution $\cZ$ of $(M\setminus \{x_o\}, J)$ are totally geodesic 
for the K\"ahler metric determined by   $ \o = d d^c_J \t$.  The induced metric on any such leaf is flat. 
\end{prop}
If $\t$  is of class $\cC^{2}$   at  $x_o$,  the K\"ahler form  $\o = d d^c_J \t$  and the  K\"ahler metric $g = d \o(\cdot, J \cdot)$ extend at  $x_o$ as    tensor fields of class $\cC^{0}$. Furthermore, by Theorem \ref{theorem53},  there exists a normalizing map $\Phi: (M, J) \longrightarrow (\bB^n, J')$ onto a manifold in normal form with   $\phi^{(0)}_{J'} = 0$. By  Lemma \ref{lemmetto},  the   2-forms $d d^c_{J'} \t_o$ and $d d^c_{\Jst} \t_o$ agree at $0$ so that $g$ is 
positive definite also in $x_o$, hence  a $\cC^0$ K\"ahler metric.
\par
\medskip
Assume now that  $\t$ is of class $\cC^k$ for some $k \geq 2$  at  $x_o$. Next theorem shows that this fact, together with the abundance of   totally geodesic flat complex curves through $x_o$,   has   important  consequences. \par
\begin{theo} \label{crucial} Let  $(M, J)$ be an $n$-dimensional  complex manifold and  $g$   a  smooth K\"ahler metric on the complement $M \setminus \{x_o\}$ of a singleton $\{x_o\}$.  Assume also that  $g$ has a positive definite $\cC^{k}$-extension at $x_o$  for  some $k \geq 2$ and that,  for  any  $0 \neq v \in T_{x_o} M$, there exists a totally geodesic, flat complex  curve $L^{v} \subset M$ with  $x_o \in L^{v}$ and $v, Jv \in T_{x_o} L^{v}$.\par
Then the curvature    $R$ of $g$ and  the  covariant derivatives   $\n^r R$,    $1 \leq r \leq k-2$,     are well-defined  at $x_o$ and they  all vanish  at  that point.  
\end{theo}
\begin{pf} Since $g$ is of class $\cC^k$, $k \geq 2$,   in any  system of coordinates around $x_o$ the  Christoffel symbols  of the Levi-Civita connection  are of class $\cC^{k-1}$ at that point and the tensors $R|_{x_o}$,  $\n^r R|_{x_o}$, $ r =  k-2$,   are well-defined. \par
\smallskip
Now,   consider a vector $0 \neq v \in T_{x_o} M$. Since the complex curve $L^{v}$ is flat and totally geodesic, for any quadruple of  vectors    $v_i \in T_x L^{v}$ at some point $x \neq x_o$, we have that  $R|_{x_o}(v_1, v_2, v_3, v_4) = 0$. 
By continuity, we get
 \beq \label{5.2} R|_{x_o}(v, Jv ,  v,  Jv ) =  0\qquad \text{for all}\ v \in T_{x_o} M\eeq
and this  implies that  $R|_{x_o} = 0$  (see e.g. \cite{KN2}, Prop. IX.7.1).  \par
\medskip
Consider now an integer $1 \leq r \leq k-2$  and assume the inductive hypothesis  $R|_{x_o} = \n R|_{x_o} = \ldots = \n^{r-1} R|_{x_o} = 0$. By Ricci formulas
(see e.g. \cite{Be}, Cor. 1.22), 
 $\n^r R|_{x_o}$  
  is  symmetric in  its first  $r$-arguments, i.e.,   
$$\n^r R|_{x_o} \in S^r(T^*_{x_o}ÊM) \otimes \left(\L^2 T^*_{x_o} M \otimes \L^2 T^*_{x_o} M\right)$$
 By $\bC$-linearity,  we may consider $\n^r R|_{x_o}$ as  a complex tensor in 
$\bigotimes^{r+4} T^{\bC*}_{x_o} M$ and, by   classical formulas of  K\"ahler geometry (see e.g. \cite{Go}, p. 166) we have that,  
for any   $w_i \in T^\bC_{x_o} M$ and  $v_j \in T^{1,0}_{x_o} M $,  
$$\n^r_{w_1 \ldots w_{r}}  R|_{x_o} ( v_1, v_2, w_{r+1}, w_{r+2}) = 0 =  \n^r_{w_1 \ldots w_r} R|_{x_o} (w_{r+1}, w_{r+2}, v_1, v_2)\ , $$
$$\n^r_{w_1 \ldots w_r} R|_{x_o} (v_1, \overline{v_2}, v_3, \overline{v_4}) =   \n^r_{w_1 \ldots w_{r}} R|_{x_o} (v_3 , \overline{v_2}, v_1, \overline{v_4}) =$$
$$ =  \n^r_{w_1 \ldots w_{r}} R|_{x_o} (v_1, \overline{v_4}, v_3, \overline{v_2})\ .$$
 From these equalities,   using Second Bianchi Identities, we also get  
 $$\n^r_{w_1 \ldots w_{r-1} v_1} R|_{x_o} (v_2, \overline{v_3}, v_4, \overline{v_5}) =   \n^r_{w_1 \ldots w_{r-1} v_2} R|_{x_o} (v_1 , \overline{v_3}, v_4, \overline{v_5}) \ ,$$
 $$ \n^r_{w_1 \ldots w_{r-1} \overline{v_1}} R|_{x_o} (v_2, \overline{v_3}, v_4, \overline{v_5}) = \n^r_{w_1 \ldots w_{r-1} \overline{v_3}} R|_{x_o} (v_2, \overline{v_1}, v_4, \overline{v_5})\ .$$
$$\n^r_{w_1 \ldots w_{r-1} \overline{v_1}} R|_{x_o} (v_2, \overline{v_3}, v_4, \overline{v_5}) = $$
$$ = \n^r_{w_1 \ldots w_{r-1} \overline{v_1}} R|_{x_o} (v_4 , \overline{v_5}, v_2, \overline{v_3}) = \n^r_{w_1 \ldots w_{r-1} \overline{v_3}} R|_{x_o} (v_2, \overline{v_1}, v_4, \overline{v_5})\ .$$
Therefore,  the value of  $\n^r R|_{x_o}$  on an ordered set of $ r +4$ vectors, $p$ of which  are  holomorphic  and  $q = r+4 -p$ are antiholomorphic,   coincides with the value   of  $\n^r R|_{x_o}$   on any other  ordered set,   given by a  permutation of   the holomorphic vectors  and a permutation of  the antiholomorphic  ones. \par
\smallskip
Consider now the following notation: given  an order set of   $r+4$ holomorphic vectors  $ w \= (w_1, \ldots, w_{r+4}) \subset  T^{1,0}_{x_o}M$, for any  $z = (z^i) \in \bC^{r+4}$ we define $w_z \= z^i w_i$.
Let also  $\cF^{(w)}(z)$ be the  homogeneous complex polynomial  of order $r+4$ in $r+4$ variables, defined by 
$$\cF^{(w)}(z) \= \n^r_{w_z  \ldots w_z} R|_{x_o}(w_z, \overline{w_z}, w_z, \overline{w_z})\ .$$
We claim that $\cF^{(w)}(z) =  0$ for any $z$. This is clear when $z$ is such that  $w_z = 0$. When  $z$ is such that $w_z \neq 0$, the vector   $w_z$ is  holomorphic  and  tangent to the flat, totally geodesic complex  curve   $L^{v}$, $v = \Re(w_z)$,  so that   $$\cF^{(w)}(z) =  \n^r_{w_z  \ldots w_z} R|_{x_o}(w_z, \overline{w_z}, w_z, \overline{w_z})  = 0\ .$$
On the other hand, by multilinearity and   symmetry properties of $\n^r R|_{x_o}$,   the  coefficient of the monomial 
$( \Pi_{i = 1}^{r+2} z^i) \overline{z^{r+3}} \overline{z^{r+4}}$  in the  expansion of  $\cF^{(w)}(z)$   
 is  equal to  
$2! (r+2)! \n^r_{w_1 \ldots w_{r}}R|_{x_o} (w_{r+1} , \overline{w_{r+3}}, w_{r+2}, \overline{w_{r+4}})$. Since $\cF^{(w)} = 0$,  this implies that 
 $\n^r_{w_1 \ldots w_{r}}R|_{x_o} (w_{r+1} , \overline{w_{r+3}}, w_{r+2}, \overline{w_{r+4}}) = 0$
  for any choice of  holomorphic vectors $w_i \in T_{x_o} ^{1,0} M$. \par
  A similar  argument  involving   the homogenous polynomial 
$$\cG^{(w)}(z) = \n^r_{w_z  \ldots w_z \overline{w_z}} R|_{x_o}(w_z, \overline{w_z}, w_z, \overline{w_z})$$
shows that   $\n^r_{w_1 \ldots w_{r-1} \overline{w_{r}}}R|_{x_o} (w_{r+1} , \overline{w_{r+3}}, w_{r+2}, \overline{w_{r+4}}) = 0$
for any choice of $w_i  \in T_{x_o} ^{1,0} M$. 
Iterating this  line of arguments, one shows that  
$\n^r R|_{x_o} = 0$. 
By induction on $r$, the theorem follows.
 \end{pf}

This  result brings directly to the next theorem. \par
\begin{theo} \label{stoll} Let  $(M,   J)$ be a manifold of circular type  with Monge-Amp\`ere exhaustion $\t$ and center $x_o$. If $\t$ is of class $\cC^{2 k}$ at $x_o$ for some integer  $k \geq 1$,  then there is a normalizing map  mapping $(M, J)$   into a manifold  in normal form $(\bB^n,   J')$ that has  $\phi^{(\ell)}  = 0$ for all $0 \leq \ell\leq k-1$. 
\end{theo}
\begin{pf} Since $\t$ is of class at least $\cC^2$ at $x_o$, by 
 Theorem \ref{theorem53}, there exists a normalizing map $\Phi: (M, J) \longrightarrow (\bB^n, J')$ for a $J'$ with   $\phi^{(0)} = 0$.  The claim is proved if we    show that, when $k \geq 2$, the almost complex structure   $J'$  satisfies also the condition    $\phi^{(\ell)} = 0$ for every $1 \leq \ell \leq k-1$.
For this,  we first remark    that, by Theorem \ref{crucial},   the K\"ahler metric $g = d d^c_{J'} \t_o(\cdot, J \cdot)$ on $(\bB^n \setminus\{0\}, J')$   admits a $\cC^{2k-2}$ extension at $0$  such that 
$\n^m R|_{0} = 0$ for any $0 \leq m\leq 2(k-2)$. Note also that such extension necessarily coincides with the standard Euclidean metric at $0$. \par
\smallskip
Pick two vectors 
$v,  w \in  S^{2n-1} \subset T_0 \bB^n$, with $w$  orthogonal to $v$  and $ \Jst v$ 
 with respect to  the standard Euclidean metric.  Note that, by means of  affine translations,  the vectors $w$ and $Jw$  can be identified with two  elements of  $T_{v} S^{2n-1}$ and such  that  $w^{1,0} \= \frac{1}{2}(w - i \Jst w)$ $ (= \frac{1}{2}(w - i  J' |_0w)) $ belongs to the CR distribution $\cH^{1,0} $ of the unit sphere $S^{2n-1}$ in $ T_0 \bB^n$.\par
Consider now the  unique stationary  disk    $f^{v}Ê: \D \to \bB^n$ of $(\bB^n, \cA_{J'})$  with 
$f^{v}(0) = 0$ and $ f^{v}_*\left(\left.\frac{\partial}{\partial x}\right|_{0}\right) = \l v$ for some  $\l > 0$, and let $\s \subset \D$ be 
a sector   around the segment $(0,1)$. Let also    $\xi = (z^j = x^j + i y^j) $ be a system of $J'$-holomorphic coordinates,  defined on a neighborhood of $f^{v}(\s)$ and  satisfying the following two conditions:  
\begin{itemize}
\item[a)]Êin such  coordinates,  the  stationary  disk    $f^{v}Ê$ is  of the form
$f^{v}(\z) = (\z, 0, \ldots, 0)$; 
\item[b)] at every $\z \in \s$, the vectors $ \left.\frac{\partial}{\partial z^j}\right|_{f^{v}(\z)}$, $2 \leq j \leq n$,  are in $\cH^\bC$.
\end{itemize}
Coordinates satisfying  (a)  and (b) can be   constructed,  by modifying   complex polar coordinates in an appropriate way.   Moreover, since the vectors $ \left.\frac{\partial}{\partial z^j}\right|_{f^{v}(\z)} $, $2 \leq j \leq n$,   are $J'$-holomorphic and  in $\cH^{\bC}$, they are of the form $ \left.\frac{\partial}{\partial z^j}\right|_{f^{v}(\z)} = E_j|_\z +\overline{ \phi(\bar E_j|_\z)}$
for some continuous family  of  $\Jst$-holomorphic vectors $E_j|_\z \in \cH^{1,0}$. With no loss of generality,  we may also  assume  that: 
\begin{itemize}
\item[c)] the vectors   $\left.\frac{\partial}{\partial z^j}\right|_{f^{v}(\z)}$, $2 \leq j \leq n$,  are prescribed in such a way  that the vectors
$Ê E_j|_\z$ have constant components with respect to  standard complex polar coordinates and
$\displaystyle\lim_{\smallmatrix t \to 0\\ t \in (0,1)\endsmallmatrix} E_2|_{t} =  w^{1,0}$.
\end{itemize}
 Let us now denote  by  $\o_o  \=  d d^c_{\Jst} \t_o$, $\o \=  d d^c_{J'} \t_o$ and  set 
 $$\o_{oi \bar j} \= \o_o\left(\frac{\partial}{\partial z^i}, \frac{\partial}{\partial \overline{z}^j}\right) \ ,\quad \o_{i \bar j} \=  \o\left(\frac{\partial}{\partial z^i}, \frac{\partial}{\partial \overline{z}^j}\right) =  \o_{oij}+ \d \o_{i\bar j}$$
 for some smooth functions $\d \o_{i \bar j}$. Note that, by   \eqref{lemmettoeq} and condition (c),  the functions $ \o_{oij}|_{f(\z)}$ are constant  and 
 \beq \label{prugne}Ê \d \o_{i\bar j}|_{f(\z)} = d d^c_{\Jst} \t_o(\overline{\phi(\overline{E_i})}, \phi(\overline E_j))|_{f^{v}(\z)}\ \text{for all}\  \z \in \s\ .\eeq
 Let also  $\G_{ij}^k$ be  the Christoffel symbols of $\n$  in the coordinates $(z^i)$  and  denote by 
 $\G_{i j|\bar k} \= \o_{\ell \bar k} \G_{ij}^\ell =  \o(\n_{\frac{\partial}{\partial z^j}}\frac{\partial}{\partial z^i}, \frac{\partial}{\partial \bar z^k})$.  From standard  facts of K\"ahler geometry, we know that 
 $\G_{i j|\bar k} = \frac{\partial \o_{j \bar k}}{\partial z^i} $. Moreover, since  $f^{v}(\D)$ is totally geodesic and flat, we also have 
 $\G_{11|\bar k}(f^{v}(\z))=  \G_{\bar 1\bar 1|k}(f^{v}(\z))  = 0$ for any $\z \in \s \subset \D$. It follows that 
\beq \label{5.22} \left.R\left( 
\frac{\partial}{\partial z^1}, \frac{\partial}{\partial \bar z^1}, \frac{\partial}{\partial z^2}, \frac{\partial}{\partial \bar z^2}\right)\right|_{f^{v}(\z)} =  \left. \frac{\partial \G_{\bar 1 \bar 2| 2}}{\partial z^1}\right|_{f^{v}(\z)} =   \left. \frac{\partial^2 \o_{2\bar 2}}{\partial z^1 \partial \bar z^1}\right|_{f^{v}(\z)} = $$
$$ =  \left. \frac{\partial^2\left( \left.\d\o_{2 \bar 2}\right|_{f^{v}(\z)}\right)}{\partial \z \partial \bar \z}\right|_{\z} =
   \left. \frac{\partial^2 \left(d d^c_{\Jst} \t_o(\overline{\phi(\overline{E_2})}, \phi(\overline E_2))|_{f^{v}(\z)}\right)}{\partial \z \partial \bar \z}\right|_{\z}  \ .\eeq
  Since   $R|_0 = 0$ and $\phi^{(0)} = 0$, this implies that   
  $$0 = \!\!\lim_{\smallmatrix t \to 0\\ t \in (0,1)\endsmallmatrix} Ê  \left. \frac{\partial^2 d d^c_{\Jst} \left(\t_o(\overline{\phi(\overline{E_2})}, \phi(\overline E_2))|_{f^{v}(\z)}\right)}{\partial \z \partial \bar \z}\right|_{t}  \!\!\! =   d d^c_{\Jst} \t_o(\overline{\phi^{(1)}_{[v]}(\overline{w^{1,0}})}, \phi^{(1)}_{[v]}( \overline{w^{1,0}}))\ .$$
 Being $dd^c_{\Jst} \t_o > 0$, it follows  that $\phi^{(1)}_{[v]}(\overline{w^{1,0}}) = 0$ for every $v, w$, hence   $\phi^{(1)} = 0$.  Due to this, the vanishing of $\phi^{(0)} = 0$ and  \eqref{prugne},  we  also have  
 \beq\label{susine} \lim_{\smallmatrix t \to 0\\ t \in (0,1)\endsmallmatrix}  \G_{1j|\bar k}(f(t)) = \lim_{\smallmatrix t \to 0\\ t \in (0,1)\endsmallmatrix}  \left.\frac{\partial( \d \o_{j\bar k})}{\partial \z}\right|_{t} = 0\qquad \text{and}$$
 $$ \lim_{\smallmatrix t \to 0\\ t \in (0,1)\endsmallmatrix}  \left.\frac{\partial \G_{1j|\bar k}}{\partial \bar z^1}\right|_{f(t)} = \lim_{\smallmatrix t \to 0\\ t \in (0,1)\endsmallmatrix}  \left.\frac{\partial( \d \o_{j\bar k})}{\partial \z \partial \bar \z}\right|_{t} = 0\qquad \text{for all}\ 1 \leq j, k \leq n\ .\eeq
From \eqref{susine},  \eqref{5.22} and the fact that $\n R|_0 = 0$, it follows that 
 \beq \label{5.21bis} 0 =\lim_{\smallmatrix t \to 0\\ t \in (0,1)\endsmallmatrix}   \n_{\frac{\partial}{\partial z^1}}\n_{\frac{\partial}{\partial \bar z^1}}
\left.R\right|_{f^{v}(t)}\left(\frac{\partial}{\partial z^1}, 
 \frac{\partial}{\partial \bar z^1}, \frac{\partial}{\partial z^2}, \frac{\partial}{\partial \bar z^2}\right)=$$
 $$ = 
  \lim_{\smallmatrix t \to 0\\ t \in (0,1)\endsmallmatrix} \left.\frac{\partial^{2}  R( 
\frac{\partial}{\partial z^1}, \frac{\partial}{\partial \bar z^1}, \frac{\partial}{\partial z^2}, \frac{\partial}{\partial \bar z^2}) }{\partial z^1 \partial \bar z^1}\right|_{f^{v}(t)} = $$
$$ = \lim_{\smallmatrix t \to 0\\ t \in (0,1)\endsmallmatrix}    \left. \frac{\partial^4 \left(d d^c_{\Jst} \t_o(\overline{\phi(\overline{E_2})}, \phi(\overline E_2))|_{f^{v}(\z)}\right)}{{(\partial \z)}^2 {(\partial \bar \z)}^2}\right|_{t} \!\!\!=  d d^c_{\Jst} \t_o(\overline{\phi^{(2)}_{[v]}(\overline{w^{1,0}})}, \phi^{(2)}_{[v]}( \overline{w^{1,0}}))\ .\eeq
Being $v$, $w$ arbitrary, this implies that $\phi^{(2)} = 0$. As before, this implies also   
\beq \lim_{\smallmatrix t \to 0\\ t \in (0,1)\endsmallmatrix}  \left.\frac{\partial^{r+s} \G_{1j|\bar k}}{(\partial z^1)^r (\partial \bar z^1)^s }\right|_{f(t)} = \lim_{\smallmatrix t \to 0\\ t \in (0,1)\endsmallmatrix}  \left.\frac{\partial( \d \o_{j\bar k})}{{(\partial \z)}^{r+1} {(\partial \bar \z)}^s}\right|_{t} = 0\eeq
 for any $r$, $s$  with  $1 \leq r + s\leq 3$ and for every indices  $1 \leq j, k \leq n$. Iterating this argument, 
 we get by induction that,   for all $0 \leq  m \leq k-2$  and for any $v$, $w$,   
$$ 0 = \lim_{\smallmatrix t \to 0\\ t \in (0,1)\endsmallmatrix} \underset{m-\text{times}}{\underbrace{\n_{\frac{\partial}{\partial z^1}} \ldots \n_{\frac{\partial}{\partial z^1}} }}
\underset{m-\text{times}}{\underbrace{\n_{\frac{\partial}{\partial \bar z^1}} \ldots \n_{\frac{\partial}{\partial \bar z^1}} }}
\left.R\right|_{f^{v}(t)}\left(\frac{\partial}{\partial z^1}, 
 \frac{\partial}{\partial \bar z^1}, \frac{\partial}{\partial z^2}, \frac{\partial}{\partial \bar z^2}\right)=  $$
 $$ = \lim_{\smallmatrix t \to 0\\ t \in (0,1)\endsmallmatrix} \left.\frac{\partial^{2m}  R( 
\frac{\partial}{\partial z^1}, \frac{\partial}{\partial \bar z^1}, \frac{\partial}{\partial z^2}, \frac{\partial}{\partial \bar z^2}) }{(\partial z^1)^m(\partial \bar z^1)^m}\right|_{f^{v}(t)} =$$
$$ =    \lim_{\smallmatrix t \to 0\\ t \in (0,1)\endsmallmatrix} \left.\frac{\partial^{2m +2} \left(\left.d d^c_{\Jst} \t_o(\overline{\phi(\overline{E_2})}, \phi(\overline E_2))\right|_{f^{v}(\z)}\right)}
{{(\partial \z)}^{m+1}  {(\partial \bar \z)}^{m+1}}  \right|_{t} =   $$
  $$ =  d d^c_{\Jst} \t_o(\overline{\phi^{(m+1)}_{[v]}(\overline{w^{1,0}})}, \phi^{(m+1)}_{[v]}( \overline{w^{1,0}})) \ .$$
This  means that  $\phi^{(\ell)}= 0$ for every $1 \leq \ell \leq k-1$.\end{pf}  
Since $J' = \Jst$ exactly when  $\phi^{(j)} = 0$ for all $j$'s,  the previous   result gives Stoll's  Theorem   as   immediate corollary. 
 \begin{cor}[\cite{St}] If $(M, J)$ is a manifold of circular type with  a Monge-Amp\`ere exhaustion $\t$ which is $\cC^\infty$ at all points,    there exists a biholomorphism $\Phi: M \longrightarrow \bB^n$ between $(M, J)$ and  the unit ball $\bB^n$ such that 
$ \t_o \circ \Phi = \t$.
\end{cor}
 \par
   \bigskip
   \section{Examples of convex domains  with \\
Monge-Amp\`ere exhaustions of prescribed regularity} \label{sei}
   \setcounter{equation}{0}
Let $(\bB^n, \cA_J)$ be a  circular manifold in normal form determined by a  deformation tensor $\phi = \sum_\ell \phi^{(\ell)}$.  
By Theorem \ref{prop41} and  Theorem \ref{stoll}, if there is an $m \geq 3$  such that 
{\it \beq \label{condition1} \phi^{(j)} = 0\ \ \text{for every}\  0 \leq j \leq 2 m -1 \ \  \text{and}\ \  \phi^{(2 m )} \neq 0\ ,\eeq}
then  
{\it  \beq \label{condition2}Ê\t_o\ \text{is of class}\ \cC^{m}\  \text{but  not of class}\ \cC^{4 m + 2}\eeq}
 in any system of $J$-holomorphic coordinates   of  $(\bB^n, \cA_J)$, 
 \par
   \medskip
A multitude of  examples  satisfying \eqref{condition1} can be  constructed.  Indeed, observe that if  \eqref{condition1}  holds,  clearly $\phi^{(0)} = 0$
and therefore the formally integrable almost complex structure $J$ of    $(\bB^n, \cA_J)$ is one of those considered by Bland and Duchamp in \cite{BD}.  By Thm. 18.2 and 18.4 of that paper, we have (for the definition of  the operator $\bar \partial_b$,  see \cite{BD, PS}):
   \begin{theo}\label{propina}  Let $f = \sum_{k \geq 1} f^{(k)}: \wt \bB^n \longrightarrow \bC$ be a smooth $\bC$-valued function, which is sum  of a series, which is uniformly converging on compacta and   with terms of the form
   $$f^{(k)}(w^a, \bar w^b, \z, \bar \z) = \z^k \wh f^{(k)}(w^a, \bar w^b) $$
   in    complex polar coordinates $(w^a, \z)$.  Let also $\phi_f = \sum_{k \geq 1}\phi^{(k)}_f$ be the formal  series   defined by  
 $$Ê\phi^{(0)}_fÊ= 0\ ,\qquad \phi^{(1)}_f = \overline\partial_b (\overline \partial_b f^{(1)})^\sharp\ ,$$
   $$\phi^{(k)}_f = \overline\partial_b \left(\overline \partial_b f^{(k)}\right)^\sharp + \frac{1}{2}Ê\sum_{\ell =1}^{k-1} h_{\cK}\left(\left[\phi^{(\ell)}, \phi^{(k - \ell)}\right]\right)\ ,\quad 1 \leq k < \infty\ ,$$
 where    $\left(\overline \partial_b f^{(j)}\right)^\sharp$ is the  unique  vector field satisfying  $ \imath_{\left(\overline \partial_b f^{(j)}\right)^\sharp} (d d^c_{\Jst}\t_o) =$  $2 i \overline \partial_b f^{(j)}$
and  
$h_{\cK}$ is the homotopy operator  defined in \cite{BD}, Lemma 16.4.\par
 Then  $\phi_f= \sum_{k \geq1}\phi^{(k)}_f$ converges uniformly on compacta and is the  deformation tensor of  a manifold  in normal form $(\bB^n, \cA_J)$. 
  Moreover, if   $f^{(j)} = 0$ for all $1 \leq j \leq m-1$,  then  $\phi_f$ is such that 
  $\phi^{(j)}_f = 0$ for all $1 \leq j \leq m-1$.
 \end{theo}
  Consider now a  function $f: \wt \bB^n \longrightarrow \bC$ of the form   
 $$f(w^a, \bar w^b, \z, \bar \z) = \z^{2m} \wh f(w^a, \bar w^b) \quad \text{for some}\quad   m \geq 3\ $$ 
 in complex polar coordinates $(w^a, \z)$.
 By  Theorem \ref{propina},  $f$ determines a  manifold in normal form $(\bB^n, \cA_J)$ satisfying \eqref{condition1} and hence with a  Monge-Amp\`ere exhaustion  of class  $\cC^{m}$ but not of class $\cC^{4m + 2}$ at the center. \par
 \medskip
 By  Remark 17.4 (iv)  in  \cite{BD}  (see also  \S 4.5 of \cite{BD1}), if  $f$ is  so that   $\left(\overline \partial_b f^{(j)}\right)^\sharp$ is sufficiently small in an appropriate Folland-Stein norm, the corresponding  manifold in normal form $(\bB^n, \cA_J)$ is  biholomorphic to  {\it a smoothly bounded, strictly linearly convex domain $D \subset \bC^n$, which is   arbitrarily close  to  $\bB^n$}, but,  of course,  not biholomorphic to the unit ball. \par
 \medskip
By formula \eqref{lilla},  this discussion  immediately gives  the following:
 \begin{cor} For any $m \geq 3$, there exists a strictly linearly convex domain $D \subset \bC^n$,  arbitrarily close to $\bB^n$,     with  squared Kobayashi distance from
a point  $x_o \in D$, which  is of  class  $\cC^{m}$ but not  $\cC^\infty$   (more precisely, not  $\cC^{4m + 2}$). \par
 \end{cor}  
\par
   \bigskip
\section{A refinement of Stoll's characterization of $\bC^n$}
\setcounter{equation}{0}
Consider a  complex manifold $M$  with   an  unbounded continuous exhaustion  $\t: M \longrightarrow [0,\infty)$  that is equal to $0$ on  a singleton $\{x_o\}$,  is of class $\cC^k$ on $M \setminus \{x_o\}$  and satisfies  conditions (ii) and (iii) 
of  the Monge-Amp\`ere exhaustions of    manifolds of circular type. Any such manifold  is called {\it unbounded manifold of circular type of class $\cC^k$} and $\t$ is called  {\it Monge-Amp\`ere exhaustion} of $M$.  
Using our setting we can now easily derive the refinement  of Stoll's characterization of $\bC^n$  (\cite{Bu, St})  advertised in the introduction.
 \begin{theo} \label{theo7.2} Let $M$ be a complex manifold  with a $\cC^4$ Monge-Amp\`ere exhaustion $\t:M \to [0,+\infty)$. Then  there exists a biholomorphism $\Phi: M\to \bC^n$ such that $\t(z)= \|\Phi(z)\|^{2}$ for all $z\in M$. In particular, $\t$ is  of  class  $\cC^\o$.
  \end{theo}

\smallskip
Before the proof we need the following remarks. For any unbounded  manifold of circular type $(M,\t)$,  each  subset $M(r) = \{\ x \in M\ : \ \t(x) < r\ \}$  is endowed with  the  exhaustion $\t_{r} \= \frac{1}{r}\t$ and can be considered as  a  bounded manifold of circular type, the only difference being  that $\t_r$ is just  $\cC^k$ away from the center and not $\cC^\infty$.  Due to this,  most properties   of    bounded  manifolds of circular type  extend to the unbounded ones. In particular, one can  directly check that
the circular representation considered in \cite{Pt1} can be  constructed for any unbounded manifold of circular type $M$ of class   $\cC^k$, $k \geq 2$, and determines 
a $\cC^{k-2}$-diffeomorphism  between  $M$ and the blow up  $\wt \bC^n$ of $\bC^n$ at the origin.  Composing such circular representation with  a diffeomorphism of $\wt \bC^n$ satisfying  (1) - (4) of \S \ref{section25},  we get a natural analogue of  Theorem \ref{existenceanduniqueness}. To state it, we first need some piece of notation.
\par
\smallskip
For any $m \geq 1$, let us  call  {\it  unbounded $\cC^m$-manifold in normal form}  any complex manifold of the form  $M = (\bC^n, \cA_J)$, given by  a non-standard atlas  $\cA_J$  of complex coordinates on $\bC^n$, which makes $M$ the  blow-down  at $0$   of a complex manifold  of the form $\wt M = (\wt  \bC^n, \cA_{\wt J})$,   where  $\cA_{\wt J}$ is a non-standard atlas satisfying the following conditions: 
\begin{itemize}  
\item[1)] any chart of $ \cA_{\wt J}$  is a map of class $\cC^m$ if expressed in terms of  the charts of the  standard differentiable manifold structure  of $\wt \bC^n$; 
\item[2)] the formally integrable almost complex structure $\wt J$ is given by a tensor field  of $\wt \bC^n$ of class $\cC^{m-1}$ in the standard manifold structure  of $\wt \bC^n$, and  satisfies   (i) - (iii) of Definition \ref{definition5.1}
\end{itemize} 
The function  $\t_o = \|\cdot \|^2$ is called  the  {\it  Monge-Amp\`ere  exhaustion of $(\bC^n, \cA_J)$}.

\begin{theo} \label{existenceanduniqueness-bis}Ê For any  unbounded manifold  of circular type $M=(M, \cA_J)$  of class $\cC^k$, $k \geq 3$, with exhaustion $\t$ and center $x_o$, there exists   a biholomorphism  $\Phi: (M, \cA_J) \longrightarrow  (\bC^n, \cA_{J'})$ onto an unbounded $\cC^{k-2}$-manifold in normal form,  mapping  the leaves of the Monge-Amp\`ere foliation of $M$ into the complex lines through the origin of $\bC^n$ and such that  $\t = \t_o \circ \Phi$ . 
 \end{theo}
 \medskip
 With this result, we are  able to prove Theorem \ref{theo7.2}.
  \par
 \medskip
\noindent{\it Proof of Theorem \ref{theo7.2}.}  By Theorem \ref{existenceanduniqueness-bis},  we may assume that  $(M, \cA_J)$  is an unbounded  $\cC^2$-manifold in normal form $(\bC^n, \cA_J)$ and  $\t = \t_o$. The claim is proved if we can show that   the tensor field $J$ on $\bC^n$  is  actually equal to $\Jst$.\par
As  for the (bounded) manifolds in normal form,   the tensor field $J$  is uniquely determined by  its deformation tensor $\phi = \phi_J$, defined in  \eqref{def1}. Since $J$ is of class $\cC^1$ on $\bC^n \setminus \{0\}$ (and the same is true for the corresponding tensor $\wt J$ on $\wt \bC^n$), also $\phi$ is of class $\cC^1$ and satisfies the same conditions 
of  the deformation tensors of bounded manifolds in normal forms. In particular, it is of the form $\phi = \sum_{k \geq 0} \phi^{(k)}$ for some tensor fields  $\phi^{(k)}$  of the form  \eqref{2.9}. Moreover,  since $\t$ is of class  $\cC^4$  at the center, the same arguments of Theorem \ref{theorem53} imply that we may assume  $\phi^{(0)} = 0$.\par
\smallskip
  Consider now  some $0 \neq v \in \bC^n$   and the  straight   line $f^{v}: \bC \longrightarrow \bC^n$, $f^{v}(\z) := v \z$.
Consider also  a  $\Jst$-holomorphic vector $E \in T_0^{1,0} \bC^n \simeq \bC^n$ of unit length  and orthogonal to $v$ and $Jv$. By construction,   the affine translation $E_y$ of $E$ at any point $y   \in  f^{v}(\bC)$ is  in  $\cH^{1,0}_y$ (here, we denote by   $\cH$ and   $\cH^{1,0} \subset \cH^\bC$ the analogues on $\bC^n$ of  the distributions $\cH$ and $\cH^{1,0}$  of   $\bB^n$, defined in  \S \ref{distrter}). Hence, by property (iii) of  deformation tensors, 
 \beq \label{stoll1} d d^c_{\Jst} \t_o(\overline{\phi(\overline E_y)}, \phi(\overline  E_y) ) < d d^c_{\Jst} \t_o( E_y, \overline  E_y ) = \|E\| = 1\ .\eeq
 Secondly, fix  a collection  $(e_i)$ of  linearly independent local generators for $\cH^{1,0}$, which are invariant under the flow of the real vector field  
$ZÊ= \Re\left(z^i \frac{\partial}{\partial z^i}\right)$ and defined on a neighborhood of $f^{v}(\bC) \setminus \{0\}$, 
and denote by $(e^j)$ the corresponding field of dual coframes in $\cH^{1,0}{}^*$.  We recall that, if the tensor field $\phi$  is written as linear combination of the tensor fields $e_i \otimes \overline{e^j}$, i.e., 
$$\phi|_{([v], \z v)} = \phi_{\overline j}^i|_{([v], \z v)}  e_i \otimes \overline{e^j}\ ,\qquad  ([v], \z v) \in \wt \bC^n \setminus \pi^{-1}(0) = \bC^n \setminus \{0\}\ ,$$
the functions  $ \phi_{\overline j}^i|_{([v], \z v)}$  are holomorphic in  $\z$ (see Prop. 4.2 (iii) in \cite{PS}).  \par
Thirdly, let $y_o \= f^{v}(1) = v$ and write   $E_{v} = E^i e_i|_{v}$.  By invariance under the flow of $Z$, we have that, at any point $y = \z v$, $\z \neq 0$,
$$e_i|_{\z v} = |\z|\cdot  e_i|_{v}\ ,\ \ \ e^i|_{\z v} = \frac{1}{|\z|} \cdot e^i|_{v}\qquad \Longrightarrow \qquad E_{\z v}Ê= \frac{1}{|\z|} E^i e_i|_{\z v}$$
 and  \eqref{stoll1} becomes
 \beq \label{stoll1bis} Ê\frac{\sum_{i,j}\left|E^i   \left.\phi^{\overline j}_i\right|_{([v], \z v)}\right|^2  \|e_i|_{  \z v}\|^2}
 {\phantom{ABC}{{|\z|}^2}^{\phantom{A^A}}}= \sum_{i,j}\left|E^i   \left.\phi^{\overline j}_i\right|_{([v], \z v)}\right|^2 \|e_i|_{v}\|^2<  1\ .\eeq
Being  $E$ arbitrary,  it follows that each  map $\z \longmapsto \phi^{\overline j}_i|_{([v], \z v)}$
 is holomorphic and uniformly bounded on $\bC$, hence constant by Liouville Theorem.  This implies  that $\phi = \phi^{(0)} = 0$, i.e. that $J = \Jst$, as we needed to show.
\square.\par
 \begin{rem} The proof shows that the theorem holds requiring  that the exhaustion $\t$  is  just of class $\cC^2$  at the center and $\cC^4$ elsewhere. In fact,  a simple modification of the above proof  shows that $(M, \cA_J)$ is mapped biholomorphically onto $\bC^n$ even if  $\t$ is just $\cC^0$ at the center.  In this case,  $\t$ is mapped into   the square $\mu^2_D$ of the Minkowski functional $\mu_D: \bC^n \to [0, + \infty)$ of a circular domain $D$. This fact  generalizes  Thm.  4.4 (i) in \cite{Pt1}. \par
 Indeed a more detailed analysis of the regularity of the circular map might provide similar results under (slightly) weaker differentiability assumptions on $\t$  as the main argument  consists of a leafwise use of Liouville Theorem.
 \end{rem}

\bigskip
\bigskip
\font\smallsmc = cmcsc8
\font\smalltt = cmtt8
\font\smallit = cmti8
\hbox{\parindent=0pt\parskip=0pt
\vbox{\baselineskip 9.5 pt \hsize=3.1truein
\obeylines
{\smallsmc
Giorgio Patrizio
Dip. Matematica e Informatica ``U. Dini''
Universit\`a di Firenze
Viale Morgani 67/a
I-50134 Firenze
ITALY
}\medskip
{\smallit E-mail}\/: {\smalltt patrizio@math.unifi.it
}
}
\hskip 0.0truecm
\vbox{\baselineskip 9.5 pt \hsize=3.7truein
\obeylines
{\smallsmc
Andrea Spiro
Scuola di Scienze e Tecnologie
Universit\`a di Camerino
Via Madonna delle Carceri
I-62032 Camerino (Macerata)
ITALY
}\medskip
{\smallit E-mail}\/: {\smalltt andrea.spiro@unicam.it}
}
}


\begin{thebibliography}{10}
\bibitem[1]{Be} A.L. Besse,
Einstein manifolds,
{\it Springer-Verlag}, 1986.
\bibitem[2]{BD} J. Bland and T. Duchamp, {\it Moduli for pointed 
convex domains\/}, Invent. Math. {\bf 104} (1991), 61--112. 
\bibitem[3]{BD1} J. Bland and T. Duchamp, {\it Contact geometry and CR structures on spheres\/}, 
in "Topics in Complex Analysis", Banach Center Pubb., vol. 31 (1995), 99--113.
\bibitem[4]{BDK} J. Bland, T. Duchamp and M. Kalka, {\it On the automorphism 
group of strictly convex domains in $\bC^n$}, Contemp. Math. {\bf 49} (1986), 19--29.
\bibitem[5]{Bu} D. Burns, {\it Curvatures of Monge-Amp\`ere foliations and 
parabolic manifolds}, Ann. of Math. {\bf 115} (1982), 349--373.
\bibitem[6]{Go} S. I. Goldberg,  Curvature and Homology, {\it Academic Press}, 1970.
\bibitem[7]{HT} C. D. Hill and M. Taylor, {\it Integrability of Rough Almost Complex Structures}, J. Geom. Anal. {\bf 13} (1) (2003), 163--172.
\bibitem[8]{KN2} S. Kobayashi and K. Nomizu, Foundations of Differential Geometry, vol. 2, {\it John Wiley \& Sons, Inc.}, 1996.
\bibitem[9]{Le}  L. Lempert, 
{\it La m\'etrique de Kobayashi et la repr\'esentation
des domaines sur la boule\/}, Bull. Soc. Math. France, 
{\bf  109} (1981), 427--474. 
\bibitem[10]{Le2}  
L. Lempert,
{\it Holomorphic invariants, normal forms and the moduli space 
of convex domains\/}, 
Ann. of Math.  {\bf 128}
(1988),  43--78. 
\bibitem[11]{Ma} B. Malgrange, ``Sur l'integrabilite des structure presque-complex'', in  Symposia Math. vol. II, 289--296, 
{\it Academic Press}, 1969. 
\bibitem[12]{NN} A. Newlander and L. Niremberg, {\it Complex analytic coordinates in almost complex manifold}, 
Ann. of Math. {\bf 65} (1957), 391--404.
\bibitem[13]{NW} A. Nijenhuis and W. Woolf, {\it Some integration problem in almost-complex manifolds}, 
Ann. of Math. {\bf 77} (1963), 424--489.
 \bibitem[14]{Pt}  G. Patrizio,  {\it Parabolic Exhaustions for Strictly Convex Domains\/}, 
Manuscripta Math. {\bf 47} (1984), 271--309. 
\bibitem[15]{Pt1}  G. Patrizio,  {\it A characterization of complex manifolds 
biholomorphic to a circular domain\/}, 
Math. Z. {\bf 189} (1985), 343--363. 
\bibitem[16]{Pt3}  G. Patrizio, {\it Disques extr\'emaux de Kobayashi et \'equation 
de Monge-Amp\`ere complex}, C. R. Acad. Sci. Paris, S\'erie I, {\bf 305} (1987), 721--724.
\bibitem[17]{PS} G. Patrizio and A. Spiro, {\it Monge-Amp\`ere equations and moduli spaces of   manifolds of circular type}, 
Adv. Math. {\bf 223} (2010), 174--197.
\bibitem[18]{Si}  N. Sibony, {\it Remarks on the Kobayashi metric}, Unpublished manuscript (1979).
\bibitem[19]{St} W. Stoll, {\it The characterization of strictly parabolic manifolds}, Ann. Scuola Norm. Sup.  Pisa, Cl. Sci. (4) {\bf 7} (1980), 87--154.
\bibitem[20]{We} S. Webster, {\it A new proof of the Newlander-Niremberg theorem}, Math. Zeit. {\bf 201} (1989), 303--316.
\end{thebibliography}
\end{document}